\documentclass{amsart}

\usepackage{subcaption}
\usepackage[utf8]{inputenc}
\usepackage{afterpage}

\usepackage{graphicx}
\usepackage{amssymb}
\usepackage{amsfonts}
\usepackage{amsmath}
\usepackage{wasysym}
\usepackage{amsmath,amscd}
\usepackage{enumerate}
\usepackage{oplotsymbl}

\newtheorem{thm}{Theorem}[section]

\newtheorem{lem}[thm]{Lemma}

\newtheorem{obs}[thm]{Observation}

\theoremstyle{definition}
\newtheorem{definition}[thm]{Definition}
\newtheorem{proposition}[thm]{Proposition}
\newtheorem{corollary}[thm]{Corollary}

\theoremstyle{remark}

\numberwithin{equation}{section}

\newcommand{\R}{\mathbb{R}}  

\begin{document}
\author{Gorka Guardiola Múzquiz}
\address{Rey Juan Carlos University, Madrid, Spain}
\email{gorka.guardiola@urjc.es}
\date{\today}

\title{$\pi_B$ in asymmetric Minkowski normed spaces}

\date{\today}

\begin{abstract}
We extend the classical results of Stanisław Gołab, on the values of pi in arbitrary
normed planes, to asymmetric norms where the unit ball has one axis of symmetry.
First, we characterize the values of $\pi_B$ for
different families of polygons as unit ball $B$.
Then we prove that $\pi_B \ge 3$, and can take all possible values and is not bounded in such spaces.
\end{abstract}

\vspace{10px}


\maketitle

\section{Introduction}

Starting with an appropriately chosen set for the unit ball, $B$,
a Minkowski functional can be defined to act as a norm. That norm is
completely characterized by the geometric properties of the set $B$. In
the constructed Banach space, the value of $\pi_B$, the half-perimeter of the unit ball $B$
measured by its norm is not $\pi=3.1415\ldots$

Classical results by Gołab state that for any norm in $\mathbb{R}^2$, $\pi_B \in [3, 4]$
and that for norms that have a quarter-turn symmetry, $\pi_B \ge \pi$.
The papers~\cite{adler, euler, keller, thompsongolab, extremal}, prove
this and other facts about $\pi_B$ for different norms in $\mathbb{R}^2$. Under certain
conditions (a unit ball with one axis of symmetry), a definition of the
measurement of the perimeter can be made for asymmetric Minkowski normed
spaces, where the unit ball is not centrally symmetric.  For these
spaces, we calculate $\pi_B$ for different sets acting as unit ball, $B$, (e.g. all the
regular polygons), and find that, in general, $\pi_B \in [3, \infty)$, taking all
possible values.

In this section and the next one, some of the results we state are
from~\cite{thompson} and~\cite{extremal} and cited without proof.
In Section~\ref{sec:asym} we introduce asymmetric norms and the conditions for which the perimeter (and thus $\pi_B$) is well defined.
In Section~\ref{sec:polyg} we characterize all the norms (including asymmetric norms) for which the unit ball $B$ is a convex regular polygon.
In Section~\ref{sec:offset} we calculate $\pi_B$ for some convex polygons, with different centers.
In Section~\ref{sec:bounds} we find the bounds for $\pi_B$ for any asymmetric norm.


Along the paper, 
$\mathbb{R}$ are the real numbers, and
$\mathbb{R}^+$ are the strictly positive real numbers $x > 0$. We use $B_1\subset B_2$ like ~\cite{thompson} to mean
what other authors write as $B_1\subseteq B_2$, that is, the subset inclusions we consider are never strict.

$B$ is the unit Ball, 
$B := \{x\in X: ||x||\le 1\}$,
its interior $B^o := \{x\in X: ||x|| < 1\}$,
and its boundary
$\partial B := \{x\in X: ||x||=1\}$.

We will denote $\pi$ the half-perimeter of the unit ball of the Euclidean norm, i.e. $\pi=3.14\ldots$ and $\pi_B$ with a subindex $B$ the half-perimeter
of the unit ball of a norm specified by $B$.

\begin{definition}\label{norm}
A norm on a real linear space $X$ (i.e. a finite dimensional real vector space) is a mapping $||\cdot||$ from X into $\mathbb{R}^+$ that satisfies
the following axioms:
\begin{enumerate}[(i)]
\item $||x|| \ge 0$ with equality if and only if $x=0$,
\item $||\alpha x|| = |\alpha| ||x||$ ($\alpha \in \mathbb{R}$),
\item $||x+y||\le ||x|| +||y|| $.
\end{enumerate}
\end{definition}

\begin{definition}\label{def:convex}
A set, $B$, is convex if it contains all the straight line segments connecting two points, $x, y \in B$.
\end{definition}
\begin{definition}\label{def:centrally}
A set is centrally symmetric if, given a point $x$ contained in the set, its antipodal point, $-x$, is also contained in the set.
\end{definition}
\begin{proposition}\label{ballnorm}
If $B$ is a centrally symmetric convex set in $X$ such that each
line through $0$ meets $B$ in a non-trivial, closed, bounded segment then $||\cdot||_B$ is a norm in $X$
where $||\cdot||_B$ comes from the next definition, given $B$ meets the properties stated above (proposition 1.1.8 in~\cite{thompson}).
\end{proposition}

\begin{definition}\label{def:funct}
The functional $||\cdot||_B$, called the \emph{Minkowski functional of $B$} is defined by
$||x||_B := inf\{\xi\in\mathbb{R}^+ : x \in \xi B\}$,
or, equivalently, given that $B$ is closed and the space is complete:
$||x||_B := \frac{1}{sup\{\xi\in\mathbb{R}^+ : \xi x \in B\}}$.
\end{definition}

\begin{thm}\label{convexball}
The unit ball in a Minkowski space given by a linear space
and a norm, ($X$, $||\cdot||_B$) is centrally symmetric and convex (theorem 1.1.6 in~\cite{thompson}).
\end{thm}

\begin{definition}\label{perimeter}
The perimeter of a polygon $\mathcal{P} := [x_0, x_1,\ldots,x_n]$, with $x_i \in \mathbb{R}^2$ with respect to a norm $||\cdot||_B$ is $\mu_B(P) := \sum_{i=1}^n ||x_i -x_{i-1}||_B$. 
\end{definition}

Given a curve $c := \{x: x = \phi(t)\, \alpha \le t \le \beta\}$
and partitions of $[\alpha, \beta]$ given by $\alpha = t_0 < t_1 < \cdots < t_n = \beta$ we may consider a family of polygons,
one for each partition
$\mathcal{P}_a= [\phi(t_0), \phi(t_1), \ldots, \phi(t_n)]$. Because of the triangle inequality, the perimeters of the polygons increase as the partition is refined.
A curve is \emph{rectifiable} if the perimeter of each of the polygons $\mathcal{P}_a$ associated to a partition of the curve is bounded from above.

\begin{definition}
The length of a rectifiable curve $c \in \mathbb{R}^2$ with an inscribed family of polygons as described above $\mathcal{P}_a$ is defined to be: $\mu_B(c) := sup_a \{\mu_B(\mathcal{P}_a)\}$.
\end{definition}

\begin{proposition}
\label{inscribed}
If a convex curve $c_1$ is inscribed inside another convex curve $c_2$, its measure is smaller,
i.e. $\mu_B(c_1) \le\mu_B(c_2)$ (Proposition 4.3.4 in~\cite{thompson}).
\end{proposition}

From~\cite{extremal} we have the following:
\begin{proposition}
\label{recticonvex}
The boundary of a convex set, $c\in\mathbb{R}^2$, is a convex path. It is a convex curve, is rectifiable and, thus, has finite measure, i.e. $\mu_B(c) \in [0, \infty)$.
\end{proposition}
By proposition~\ref{recticonvex}, the length of the unit ball boundary $\partial B$, which is a convex path
by~\ref{convexball}, is finite, i.e. $\mu_B(\partial B) \ge k$, $k \in \mathbb{R}$. Also $k > 0$ because of the other conditions for the set $B$ in Proposition~\ref{ballnorm}.
\section{Asymmetric Minkowski normed spaces\label{sec:asym}}

We are interested in \emph{asymmetric} norms. That means we relax axiom (ii) in \ref{norm}.
\begin{definition}
\label{def:asymnorm}
An asymmetric norm $||x||$ has the following properties.
\begin{enumerate}[(i)]
\item $||x|| \ge 0$ with equality if and only if $x=x_0$,
\item $||\alpha x|| = \alpha ||x||$ ($\alpha \in \mathbb{R^+}$),
\item $||x+y||\le ||x|| +||y|| $.
\end{enumerate}
\end{definition}
All the properties are the same of a norm, except $||x||$ and $||-x||$ may be different. The origin is not $0$ but an
arbitrary $x_0 \in B^o$.

\begin{proposition}\label{asymballnorm}
If $B$ is a convex set in $X$ such that each
line through $x_0$ meets $B$ in a non-trivial, closed, bounded segment (i.e. a segment
of positive bounded length with edges different than the origin), then $||\cdot||_B$ is a asymmetric norm in $X$
where $||\cdot||_B$ comes from the same Minkowski functional defined by \ref{def:funct}.
\end{proposition}

Note that we have dropped the condition of central symmetry in $B$. We will prove Proposition~\ref{asymballnorm} for the
appropiate off-center Minkowski functional next.
In the original proposition \ref{ballnorm}, the center from which to calculate the norm was defined implicitly. As the set was
centrally symmetric around $0$, we could use that center to define the Minkowski functional. In the asymmetric case, though, there is
more freedom in the choice of center. Any center $x_0 \in B^o$ is equally valid:

\begin{definition}\label{def:minkcent}
The functional $||\cdot||_{B, x_0}$, called the \emph{offset Minkowski functional of $B$ at $x_0$} is defined by
$||x||_{B, x_0} := inf\{\xi\in\mathbb{R}^+ : x \in \xi (B - x_0)\}$
or, equivalently,
$||x||_{B, x_0} := \frac{1}{ sup\{\xi\in\mathbb{R}^+ : \xi x \in (B - x_0)\}}$,
where $B - x_0$ means: ``the set composed of all elements $x-x_0$ where
$x \in B$''.
\end{definition}

\begin{thm}
An asymmetric norm defined by its Minkowski functional in~\ref{def:minkcent} meets the properties stated in Definition \ref{def:asymnorm}.
\begin{proof}
Linearity $(i)$ comes trivially from the definition of the functional. Non-negativity (ii)
comes from the fact that $x_0 \in B^o$ and that  each line through $x_0$ meets
B in a non-trivial closed bounded segment.
We will use $||x||_{B, x_0}$ or $||x||_B$ when the choice of center is obvious or does not matter.

The triangle inequality (iii) is met because of the convexity of the unit ball (we can just reuse
the classic proof):
Given the definition in~\ref{def:funct}, $||x||_B := \frac{1}{sup\{\xi\in\mathbb{R}^+ : \xi x \in B\}}$,
a vector $x$ is inside the unit ball, i.e. $x \in B$, if and only if $||x||_B \le 1$. 
For a convex set $B$, given two points, $x$ and $y$, any vector $z=(1-t)x + ty$, where $t\in[0,1]$ then
$z \in B$ (this is the standard definition of a convex set).
This means that
$\frac{x+y}{||x||_B+||y||_B} =  \left(\frac{x}{||x||_B}\right)\frac{||x||_B}{||x||_B+||y||_B}  + \left(\frac{y}{||y||_B}\right)\frac{||y||_B}{||y||_B+||y||_B} =
t\frac{x}{||x||_B}  + (1-t)\frac{y}{||y||_B}$,
where $t=\frac{||x||_B}{||x||_B+||y||_B}$.
Clearly, $\frac{x+y}{||x||_B+||y||_B} \in B$, so $\frac{||x+y||_B}{||x||_B+||y||_B} \le 1$ and
 $||x+y||_B \le ||x||_B+||y||_B.$
\end{proof}
\end{thm}

When some families of sets are used as unit ball (like regular polygons), there is an obvious choice of center.
In other cases (or even for regular polygons), different centers can be defined. We will try to tackle this
question with as much generality as possible (in $\mathbb{R}^2$).

Another important obstacle when dealing with asymmetric norms in $\mathbb{R}^2$ is that the definition
of perimeter in Definition~\ref{perimeter} is not unique and depends on the choice
of direction for each of the sides.  In general 
$\sum_{i=1}^n ||x_i -x_{i-1}||_B \neq \sum_{i=1}^n ||x_{i-1} -x_i||_B $.
Even worse than two possible measurements of the perimeter of a polygon,
(clockwise or counterclockwise) we can obtain, potentially, $2^n$ different measurements
for the perimeter: two values per side, one per orientation.
Clearly, some additional conditions must be met for the half-perimeter (i.e. $\pi_B$) to be unique.

First, we give some definitions and then we will state the conditions:

\begin{definition}\label{leftperimeter}
The counterclockwise perimeter of a polygon $\mathcal{P} := [x_0, x_1,\ldots,x_n]$ , $x_i\in\mathbb{R}^2$ with respect to a possibly asymmetric norm $||\cdot||_B$ for $\mathbb{R}^2$ is
$\overset{\circlearrowleft}{\mu}_B(P) := \sum_{i=1}^n ||x_i -x_{i-1}||_B$
where the enumeration of vertices of the polygon is ordered counterclockwise.
\end{definition}
\begin{definition}\label{rightperimeter}
With the same enumeration, but different order, the clockwise perimeter is:
$\overset{\circlearrowright}{\mu}_B(P) := \sum_{i=1}^n ||x_{i-1} -x_i||_B$.
where the enumeration of vertices of the polygon is ordered clockwise.
\end{definition}

\begin{definition}\label{maxperimeter}
The maximum perimeter of a polygon $\mathcal{P} := [x_0, x_1,\ldots,x_n]$, $x_i\in\mathbb{R}^2$ with respect to a possibly asymmetric norm $||\cdot||_B$
for $\mathbb{R}^2$ is $\overset{max}{\mu}_B(P) := \sum_{i=1}^n max\{||x_i -x_{i-1}||_B, ||x_{i-1} -x_i||_B\}$.
\end{definition}

The maximum perimeter is equivalent to taking the perimeter with respect to a (symmetric)
norm obtained from the symmetrization of $B$. 
$\overset{max}{\mu}_B(P) =\mu_{B_{si}}(P) $ where $B_{si} = B \cap -B$ and $-B$ denotes the set $-B := \{-x : x \in B\}$.
Note that a smaller $B$ radius means a bigger measurement of the corresponding vector.
$B_{si}$ is the intersection of two convex sets (which is itself convex), so it meets all the conditions to be the unit ball of an asymmetric norm
as stated in Proposition~\ref{asymballnorm}. $B_{si}$, is also centrally symmetric,
from its definition, so $||\cdot||_{B_{si}}$ is a norm.

\begin{definition}\label{minperimeter}
The minimum perimeter of a polygon $\mathcal{P} := [x_0, x_1,\ldots,x_n]$ with respect to a possibly asymmetric norm $||\cdot||_B$ for $\mathbb{R}^2$ is $\overset{min}{\mu}_B(P) := \sum_{i=1}^n min\{||x_i -x_{i-1}||_B, ||x_{i-1} -x_i||_B\}$.
\end{definition}

Unfortunately the union of convex sets is not necessarily convex, so
\mbox{$B_{su} = B\, \cup\, -B$} is not convex so $||\cdot||_{B_{su}}$ may not be a
norm. An example of this is a triangle. The union of the triangle and
its reflection is a star polygon (an hexagram), non-convex.  The convex hull of a set $A$, $conv(A)$,
is the minimum
convex set containing it. The smallest symmetric ball containing $B_{su}$
is its convex hull $conv(B_{su})$.
Note that if the unit ball of a norm (asymmetric or not) is contained inside another $B_1 \subset B_2$, then 
$||v||_{B_2} \le ||v||_{B_1} $ from the definition of the Minkowski functional and the triangle inequality.

From all these definitions, it is evident that:
\begin{equation}
\mu_{conv(B_{su})}(P) \le \overset{min}{\mu}_B(P) \le \overset{\circlearrowleft}{\mu}_B(P) \le \overset{max}{\mu}_B(P) =\mu_{B_{si}}(P)
\end{equation}
 and 
\begin{equation}
\mu_{conv(B_{su})}(P) \le \overset{min}{\mu}_B(P) \le \overset{\circlearrowright}{\mu}_B(P) \le \overset{max}{\mu}_B(P) =\mu_{B_{si}}(P).
\end{equation}

For a symmetric norm
$\overset{max}{\mu}_B(P) = \overset{min}{\mu}_B(P)$
so $\overset{\circlearrowleft}{\mu}_B(P) = \overset{\circlearrowright}{\mu}_B(P)$.
For a strictly asymmetric norm (an asymmetric
norm with a not centrally symmetric unit ball\footnote{Remember all norms are, by
definition an asymmetric norms which have to meet a subset of the conditions
met by a norm.})
$\overset{min}{\mu}_B(P) \leq \overset{max}{\mu}_B(P)$.
Note that, while  $\overset{min}{\mu}_B(P) = \overset{max}{\mu}_B(P)$ implies $\overset{\circlearrowleft}{\mu}_B(P) = \overset{\circlearrowright}{\mu}_B(P)$,
the converse is not true. We want to impose the minimum conditions necessary, so that
$\overset{\circlearrowleft}{\mu}_B(P) = \overset{\circlearrowright}{\mu}_B(P)$.

\begin{figure}[htb]
	\centering
	\includegraphics[width=0.5\textwidth]{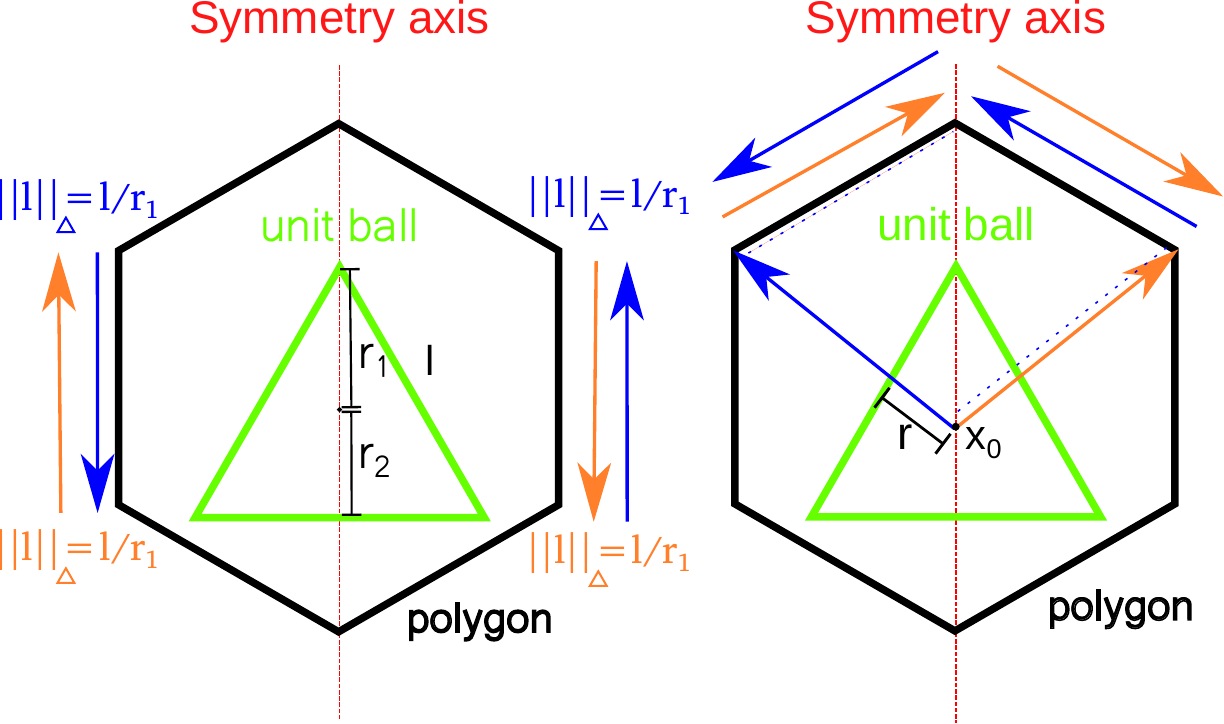}
	\caption{Ball and polygon with an axis of symmetry\label{fig:vert}}
\end{figure}

Under what conditions is $\overset{\circlearrowleft}{\mu}_B(P) =
\overset{\circlearrowright}{\mu}_B(P)$?  Just one condition is sufficient:
a single axis of symmetry shared by the unit ball and the polygon whose
perimeter is being measured.  Under this condition, each side (or fragment
of a side if the axis of symmetry cuts a side in two), can be paired with
its symmetric, so each side measured clockwise will be transposed for the
same side when measured counterclockwise but with its symmetric radius of the
same length as depicted in Figure~\ref{fig:vert}.

If we want to measure a convex path $C$,
and it shares the same axis of symmetry, we can rectify it with polygon
refinements,  $P_i$, with a growing number of sides. The refinements grow as
more points are added because of the triangle inequality so $\mu_B(P_i)
\ge \mu_B(P_j)$ when $i \ge j$, but as it is bounded because it is
a convex path, the length of the refinements converges. We will denote  the limit of the measure
$\mu_B(C) = sup \,\mu_B(P_i)= sup \,\overset{\circlearrowleft}\mu_B(P_i)=
sup \,\overset{\circlearrowright}\mu_B(P_i)$.

\subsection{Radon norms}

From~\cite{extremal} we introduce the following definitions:

\begin{definition}\label{radon}
Let $X$ be a norm on $\R^n$ and consider two nonzero vectors \\ $x, y \in \R^n$.
We will say $y$ is \emph{Birkhoff orthogonal} to $x$, written $y \vdash_X x$, if \\ $||x+ty||_X \ge ||x||_X$
for all $t\in \R$. This is equivalent to stating that the line $\{x+ty:t\in\R \}$ is ``tangent''
to the ball of radius $||x||_X$ meaning it intersects the  $\partial B$. but not $B^o$.
\end{definition}
For a Euclidean norm we have $y \vdash_X x$, if and only if $\langle x,y\rangle = 0 $, so $\vdash_X$ is
symmetric. This symmetry does not hold in general: $y \vdash_X x \nRightarrow x \vdash_X y.$

\begin{definition}
A norm $X\in \R^2$ is said to be~\emph{Radon} if the relation $\vdash_X$ is symmetric\footnote{Note that Radon norms are special in two dimensions. In three or more dimensions, a norm is Radon if and only if it is
Euclidean (see~\cite{extremal}), whereas in two dimensions, they can be non-Euclidean.}.
\end{definition}

This is Theorem 4.7 in~\cite{extremal} (see~\cite{thompsongolab},\cite{mink1}, \cite{mink2} and~\cite{martini2006antinorms}
for more details):

\begin{definition}\label{radonpi}
Lets call the set of all Radon norms  with a unit ball $B$, $\mathcal{R}$. Then $\pi_B \in [3,\pi]$, with
$\pi_B =\pi$, if and only if $B$ is an ellipse (i.e. it belongs to the Euclidean norm under linear isomorphisms, see~\cite{extremal}).
\end{definition}

\section{Regular convex polygons\label{sec:polyg}}

All regular convex polygons have at least two axes of
symmetry~\cite{coxeter}, so they meet the conditions stated above
(when acting both as unit ball and as polygon being measured,
i.e. to find $\pi_P$) so that
$\overset{\circlearrowleft}{\mu}_P(P) = \overset{\circlearrowright}{\mu}_P(P)$,
which we will denote simply
${\mu}_P(P)$, so $\pi_{P} = \frac{\mu_P(P)}{2}$.  Note that polygons
with an even number of sides have central symmetry~\cite{coxeter}
and as such the Minkowski functional is a norm. Polygons with an odd
number of sides are not centrally symmetric so they
induce a Minkowski functional which is an asymmetric
norm.  For any regular convex
polygon we can use as center the point where all the axes of
symmetry meet.  Even more, each side is measured, because of the
rotational symmetries~\cite{coxeter} both of the unit ball,
and the polygon being measured (both are the same) with
equal radius, so their measures are equal.  The half-perimeter of
a regular convex polygon is  then
$\pi_{B} = \frac{\mu_{B}(P)}{2} = \frac{\sum_1^n||l_i||_{B}}{2} = \frac{n||l_1||_{B}}{2}$.

We will denote it as $\pi_n$ or, for some
special cases like, when $B$ is a regular triangle, as $\pi_{\triangle}$, so we will use indistinctively $\pi_{B}$, 
$\pi_3$ or $\pi_{\triangle}$.

To illustrate the procedure, let's start with the regular triangle, 
$\pi_{\triangle} = \frac{\mu_\triangle(\triangle)}{2}$.

\begin{figure}[htb]
\begin{subfigure}{0.5\textwidth}
	\centering
	\includegraphics[width=0.5\textwidth]{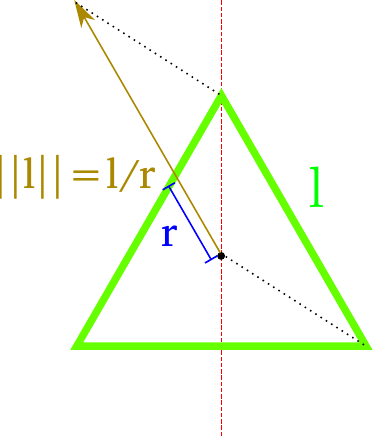}
\end{subfigure}%
\begin{subfigure}{0.5\textwidth}
	\centering
	\includegraphics[width=0.6\textwidth]{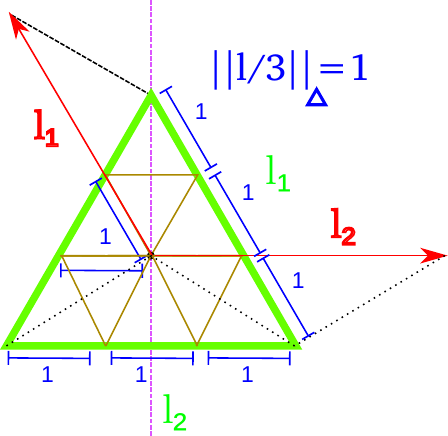}
\end{subfigure}
	\caption{Triangle as unit ball\label{fig:triang}}
\end{figure}

We will measure the perimeter anti-clockwise. We apply definition~\eqref{def:funct} to calculate $||l_n||_{\triangle}$.
For each side, we translate it to the origin to find the radius which applies to it ($r$ as can be seen
in Figure~\ref{fig:triang}). Because the triangle is the unit ball $r=1$. In Figure~\ref{fig:triang} we
measure two sides (as an example, as stated before, with one is enough, the rest are the same).
From the geometric construction
in Figure~\ref{fig:triang}, $l_1 = 3$, so $\pi_{\triangle} = \frac{\mu_\triangle(\triangle)}{2} = 3\frac{||l_1||_{\triangle}}{2} = \frac{9}{2} = 4.5$.
Note that scaling the unit ball does not change its measure, because if the length of the sides is scaled by $a$, so
is the radius by which they are being measured, by similarity, and $||l||_b = \frac{||al||}{||ar||} = \frac{||l||}{||r||}$, so $||l_b||$ is unchanged,
where $||x||$ without subindex denotes the regular Euclidean norm.\footnote{This can be generalized to a linear isomorphism, when
the polygon defines a (symmetric) norm see~\cite{extremal}. The same can be done for linear isomorphisms with positive eigenvalues for asymmetric norms.}

Note that this value is outside of Gołab's theorem range~\cite{golab}  for symmetric two dimensional norms
(i.e. $\pi_{B_{sym}} \in [3, 4]$).

When measuring the perimeter of a regular convex polygon (with itself as unit ball),
the measurements defined in~\ref{leftperimeter}, \ref{rightperimeter}, \ref{maxperimeter} and \ref{minperimeter} are the same so $max\{||l_1||_{P}, ||-l_1||_{P}\} = ||l_1||_{P} = ||-l_1||_{P}$ and $min\{||l_1||_{P}, ||-l_1||_{P}\} = ||l_1||_{P} = ||-l_1||_{P}$, and
$\overset{min}{\mu}_P(P) = \overset{\circlearrowleft}{\mu}_P(P) = \overset{\circlearrowright}{\mu}_P(P) =\overset{max}{\mu}_P(P)$:
any of these measurements would give the same value of $\pi_P$.

\begin{figure}[htb]
	\centering
	\includegraphics[width=0.5\textwidth]{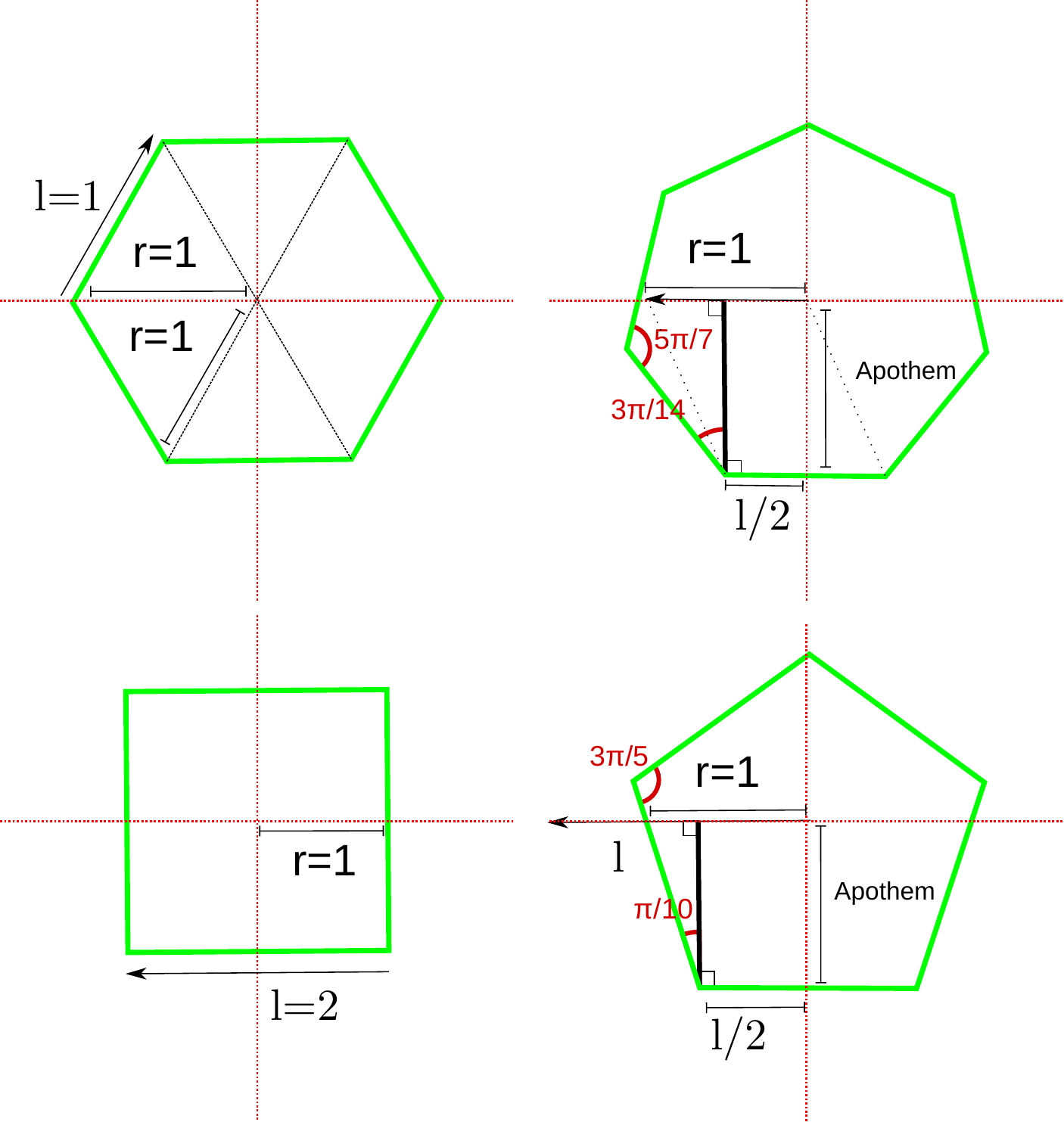}
	\caption{Constructions to calculate the measure of a side\label{fig:polyg}}
\end{figure}

Just applying some trigonometry and definition~\eqref{def:funct}, the half-perimeter of any regular
polygon can be calculated\footnote{Note that the apothem (or radius of the
inscribed circle) is: $a_n=\frac{l}{2\tan(\frac{\pi}{n})}$}, as seen in Figure~\ref{fig:polyg}.

\begin{table}[hbt]
\begin{tabular}{lll}
Number of sides &  $\pi_n$ & $\pi_n$ \\
\hline
3 & $9/2$ & $9/2$\\
4 & 4 & 4 \\
5 & 3.454915028125262879488532914085904706 & $\frac{5(5-\sqrt{5})}{4}$\\
6 &  3& 3 \\
7 & 3.286503763716470257372386327536920665 & $28\sin^2(\frac{\pi}{7})\sin(\frac{3 \pi}{14})$\\
8 & 3.313708498984760390413509793677584628 & $8\sqrt{2}-8$\\
9 & 3.225966377139231493977618069073666967 & $2+2\cos(2\frac{\pi}{9})$\\
10 & 3.090169943749474241022934171828190589 & $\frac{5+5\sqrt{5}}{2}$\\
 & & \\
\end{tabular}
	\caption{Values of $\pi_n$ vs. number of sides\label{tab:pis}}
\end{table}

\begin{figure}[htb]
	\centering
	\includegraphics[width=0.7\textwidth]{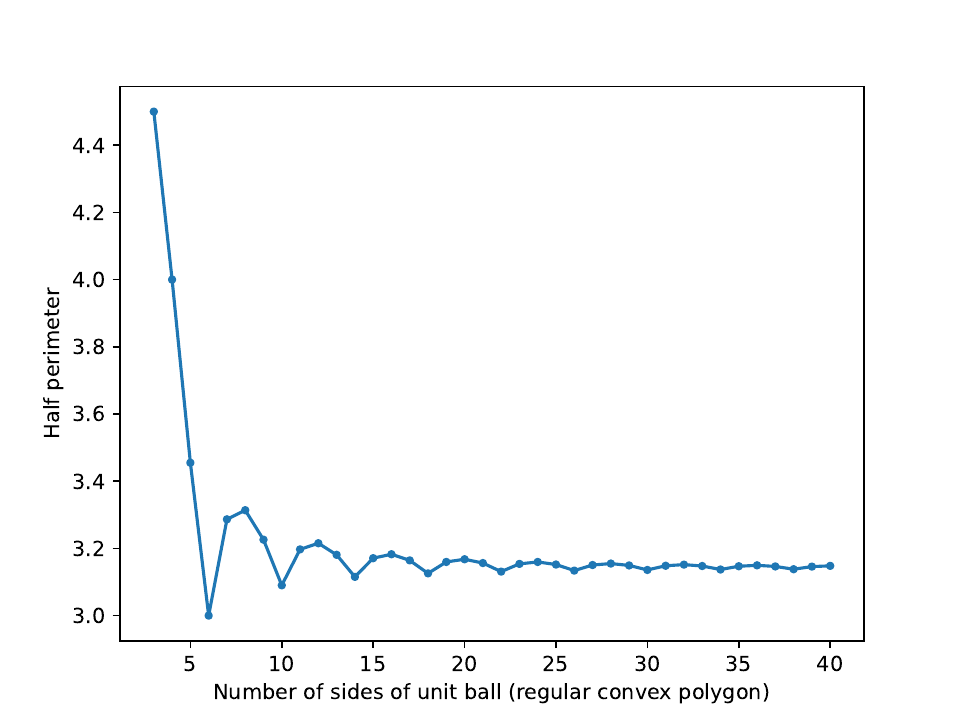}
	\caption{Values of $\pi_{B_n}$ vs. number of sides\label{fig:pis}}
\end{figure}

In Figure~\ref{fig:pis} and Table~\ref{tab:pis} we present the values of the half-perimeter for different convex
regular polygons.
A general formula can be calculated for \emph{all} the convex regular polygons. In~\cite{circlenumbers}
a formula for the Minkowski functional of a polygon is found.
\begin{equation}\label{eq:max}
||r||_{B_n}  = max_i\{\langle r, a_{(m,n)}\rangle,  i \in {1,\ldots, m}\}, 
\end{equation}
 where $\langle r, a_{(m,n)}\rangle$ denotes the inner product of the vector $r = (x, y)$
and $a_{(m,n)} = \frac{1}{\cos\left(\frac{\pi}{n}\right)}\left(\cos\left((2m-1)\frac{\pi}{n}\right),\sin\left((2m-1)\frac{\pi}{n}\right)\right)$.

The side of a polygon is the vector $l = \left(\cos(2\pi/n)-1,\sin(2\pi/n)\right)$,
so
\begin{equation*}
||l||_n = max_m\left\{\frac{\left(\cos\left(\frac{2\pi}{n}\right)-1\right)\cos\left((2m-1)\frac{\pi}{n}\right)+\sin\left(\frac{2\pi}{n}\right)\sin\left((2m-1)\frac{\pi}{n})\right)}{\cos\left(\frac{\pi}{n}\right)},  m \in {1,\ldots, n}\right\}, 
\end{equation*}
and its half-perimeter is, because of rotational symmetry,
$\pi_n=\frac{n||l||_n}{2}$,
so
\begin{equation*}
\pi_n= n\frac{max_m\left\{\left(\left(\cos\left(\frac{2\pi}{n}\right)-1\right)\cos\left((2m-1)\frac{\pi}{n}\right)\right)+\sin\left(\frac{2\pi}{n}\right)\sin\left((2m-1)\frac{\pi}{n})\right),  m \in {1,\ldots, n}\right\}}{2\cos\left(\frac{\pi}{n}\right)}.
\end{equation*}
Applying product-to-sum trigonometric identities:
\begin{equation}\label{eq:pimax}
\pi_n=\frac{n}{2\cos\left(\frac{\pi}{n}\right)}max_m\left\{\left(  \cos\left( \frac{3\pi}{n} - \frac{2m\pi}{n} \right)- \cos\left( \frac{\pi}{n} - \frac{2m\pi}{n} \right)\right),  m \in {1,\ldots, n}\right\}.
\end{equation}

\begin{figure}[htb]
	\centering
	\includegraphics[width=0.4\textwidth]{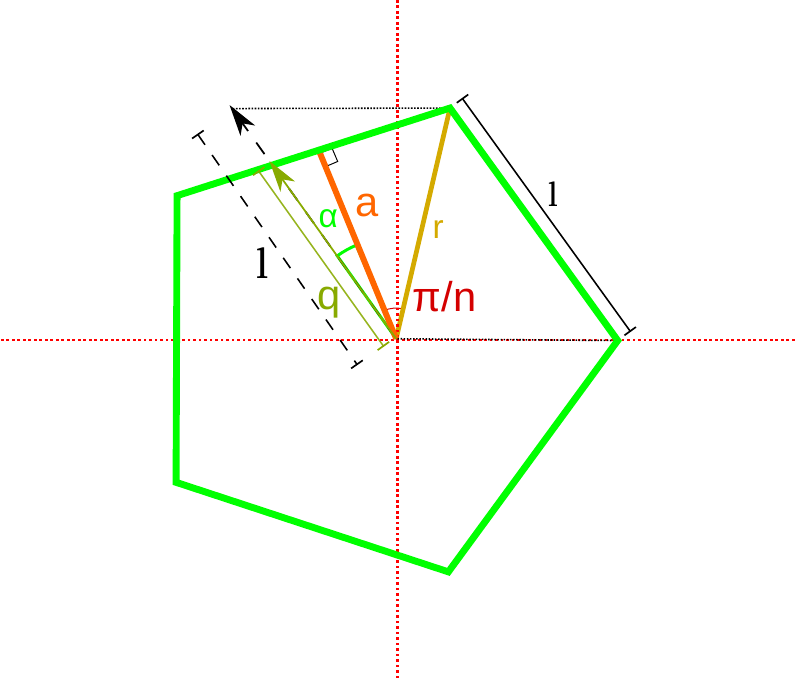}
	\caption{General construction to measure the side\label{fig:hk}}
\end{figure}

\begin{figure}[htb]
	\centering
	\includegraphics[width=0.4\textwidth]{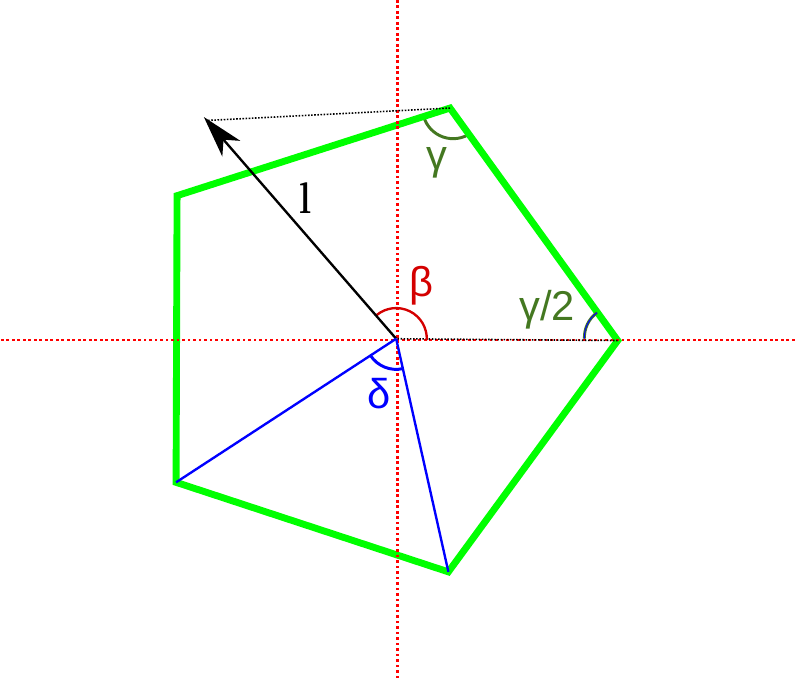}
	\caption{Sides subtended by $\beta$ \label{fig:ceil}}
\end{figure}

We can simplify this expression further. First we need to show where the expression comes from.
Figure~\ref{fig:hk} depicts the elements necessary for an example pentagon.
If, in that figure $r=1$, then $||l||_n = \frac{||l||}{||q||}$ and $||q|| = \frac{\cos(\pi/n)}{\cos(\alpha)}$.
So finally, $||l||_n = \frac{||l||\cos(\alpha)}{\cos(\pi/n)}$. If we take $u$ to be the unit vector $u= \frac{a}{||a||}$, then
$\langle u, l\rangle = ||l||cos(\alpha)$, but $u$ for the side of index $m$ is:
$u_{m,n} =  \left(\cos\left(\frac{m2\pi}{n}-\frac{\pi}{n}\right),\sin\left(\frac{m2\pi}{n}-\frac{\pi}{n}\right)\right)$.
so  $||l||_n = \frac{||l||\cos(\alpha)}{\cos(\pi/n)} =  \frac{\langle u_{m,n}, l\rangle }{\cos(\pi/n)}$, but $a_{m,n}=\frac{u_{m,n}}{\cos(\pi/n)}$, so we recover~\eqref{eq:max}.
The maximum in the equation is there to find which side is intersected by the vector (which when calculating the half-perimeter is also a side of the polygon) being measured by the functional.
We can calculate more precisely which side is intersected.
It depends on the number of sides subtended by the angle $\beta$ in Figure~\ref{fig:ceil}.
First we find the values of the different angles on the figure. The inside angle
of a regular convex polygon of $n$ sides is:
$\gamma = \frac{(n-2)\pi}{n}$.
Also, the angle $\delta$ is~\cite{coxeter}:
$\delta = \frac{2\pi}{n}$.
Finally, the number of subtended sides $m = \Big\lceil{\frac{\beta}{\delta}}\Big\rceil = \Big\lceil{\frac{\pi - \gamma/2}{\delta}}\Big\rceil =\
 \Big\lceil{\frac{\pi - \frac{(n-2)\pi}{2n}}{\frac{2\pi}{n}}}\Big\rceil =\Big\lceil{\frac{n+2}{4}}\Big\rceil\
=\Big\lfloor{\frac{n+5}{4}}\Big\rfloor$.

Substituting in~\eqref{eq:pimax} and simplifying:
\begin{equation}\label{eq:polygdef}
\pi_n=\frac{n}{2\cos\left(\frac{\pi}{n}\right)}  \left(\cos\left( \frac{\pi}{n} \left(-1+ 2\Big\lfloor{\frac{n+1}{4}}\Big\rfloor\right) \right) -\
\cos\left( \frac{\pi}{n} \left(1+2\Big\lfloor{\frac{n+1}{4}}\Big\rfloor\right) \right) \right).
\end{equation}

The limit when the number of sides $n$ tends to $\infty$ is $\pi$ as expected.
To obtain it, note that, in the limit, we can ignore the floor function as it is asymptotically negligible:
$\lim_{n \to \infty} \pi_n = \lim_{n \to \infty} \frac{n}{2\cos\left(\frac{\pi}{n+1}\right)} \left(\cos\left( \frac{\pi}{n+1} \left(-1+n/2\right) \right) -\
\cos\left( \frac{\pi}{n} \left(1+n/2\right) \right) \right) =$
$\lim_{n \to \infty} \frac{n}{2} \left(\sin\left(\frac{\pi}{n} \right) +\
\sin\left( \frac{\pi}{n} \right) \right) $,
so
$= \lim_{n \to \infty} \frac{n}{2} \left(\frac{\pi}{n} +\
 \frac{\pi}{n} \right) = \pi $.

Note also that $\pi_n$ as a function of $n$ oscillates with period 4 as seen in Figure~\ref{fig:pis}.
The function behaves as a convex function until an extra side is subtended by the
angle $\beta$. Each time this happens, it reaches a minimum and bounces back.


\begin{table}[hbt]
\begin{tabular}{ll}
Number of sides &  $\pi_n$ \\
\hline
$[3,6]$& $\frac{n}{2\cos\left(\frac{\pi}{n}\right)}  \left(\cos\left( \frac{\pi}{n}\right) -\cos\left( \frac{3\pi}{n}\right) \right) $\\
$[7,10]$ & $\frac{n}{2\cos\left(\frac{\pi}{n}\right)}  \left(\cos\left( \frac{3\pi}{n}\right) -\cos\left( \frac{5\pi}{n}\right) \right) $\\
$[11,14]$& $\frac{n}{2\cos\left(\frac{\pi}{n}\right)}  \left(\cos\left( \frac{5\pi}{n}\right) -\cos\left( \frac{7\pi}{n}\right) \right) $\\
$[15,18]$& $\frac{n}{2\cos\left(\frac{\pi}{n}\right)}  \left(\cos\left( \frac{7\pi}{n}\right) -\cos\left( \frac{9\pi}{n}\right) \right) $\\
\end{tabular}
	\caption{Values of $\pi_n$ vs. number of sides\label{tab:pisexpr}}
\end{table}

When $n=6$ the radius is equal to the side as can be seen
in Figure~\ref{fig:polyg} and the subtended number of sides
is a round number. The cosine in the right side of the expression in table~\ref{tab:pisexpr},
$\cos\left( \frac{3\pi}{n}\right) $ is zero\footnote{Something similar happens when $n=10, 14\ldots$ This is the reason these polygons
induce Radon measures as we will see later.}.

For $n$ multiple of $4$, $\Big\lfloor{\frac{n+1}{4}}\Big\rfloor=\frac{n}{4}$, and equation~\eqref{eq:polygdef} can be simplified
$\pi_n=\frac{n}{2\cos\left(\frac{\pi}{n}\right)}  \left(\cos\left( \frac{\pi}{n} \left(-1+ 2\frac{n}{4}\right) \right) -\
\cos\left( \frac{\pi}{n} \left(1+2\frac{n}{4}\right) \right) \right)$,
so:
\begin{equation}\label{eq:tanpi}
\pi_n=n \tan\left(\frac{\pi}{n}\right)
\end{equation}
For $n$ an exponent of $2$, $n=2^m$ the paper \cite{nested} gives a formula for $\sin\left(\frac{\pi}{2^m}\right) = \frac{1}{2}\sqrt{2 - \underbrace{\sqrt{2 + \sqrt{2}}\cdots}_\text{m-2 times}}$ 
and $\cos\left(\frac{\pi}{2^m}\right) = \frac{1}{2}\sqrt{2 + \underbrace{\sqrt{2 + \sqrt{2}}\cdots}_\text{m-2 times}}$, so, for $m \ge 2$
so $\pi_{2^m} = 2^m \frac{\sqrt{2 - \overbrace{\sqrt{2 + \sqrt{2}}\cdots}^\text{m-2 times}}}{\sqrt{2 + \underbrace{\sqrt{2 + \sqrt{2}}\cdots}_\text{m-2 times}}}$.
From this we can obtain a Viète-like formula:
\begin{equation}
\pi = \lim_{m \to \infty} 2^m \frac{\sqrt{2 - \overbrace{\sqrt{2 + \sqrt{2}}\cdots}^\text{m-2 times}}}{\sqrt{2 + \underbrace{\sqrt{2 + \sqrt{2}}\cdots}_\text{m-2 times}}}.
\end{equation}
\subsection{Simpler expression}

We start with equation~\eqref{eq:polygdef} and we write $n=4m+k$, where $m \in \{1,2,\ldots\}$
and $k\in \{0,1, 2, 3\}$:

\begin{equation*}
\pi_n=\frac{n}{2\cos\left(\frac{\pi}{n}\right)}  \left(\cos\left( \frac{\pi}{n} \left(-1+ 2\Big\lfloor{\frac{4m+k+1}{4}}\Big\rfloor\right) \right) -\
\cos\left( \frac{\pi}{n} \left(1+2\Big\lfloor{\frac{4m+k+1}{4}}\Big\rfloor\right) \right) \right) = 
\end{equation*}
\begin{equation*}
\frac{n}{2\cos\left(\frac{\pi}{n}\right)}  \left(\cos\left( \frac{\pi}{n} \left(-1+ 2m+2\Big\lfloor{\frac{k+1}{4}}\Big\rfloor\right) \right) -\
\cos\left( \frac{\pi}{n} \left(1+2m+2\Big\lfloor{\frac{k+1}{4}}\Big\rfloor\right) \right) \right).
\end{equation*}
Applying trigonometric angle sum identities:
\begin{equation*}
\pi_n=\frac{n}{2\cos\left(\frac{\pi}{n}\right)} 2\sin\Big(\frac{2\pi}{n}\Big(m+\Big\lfloor{\frac{k+1}{4}}\Big\rfloor\Big)\Big)sin\Big(\frac{\pi}{n}\Big) = n\tan\Big(\frac{\pi}{n}\Big)\sin\Big(\frac{2\pi}{n}\Big(m+\Big\lfloor{\frac{k+1}{4}}\Big\rfloor\Big)\Big).
\end{equation*}
When $n=4m$, $k=0$, the right side equals $\sin(\pi/2) = 1$ and we recover equation~\eqref{eq:tanpi}.
We have three other cases, $k=1$, $k=2$ and $k=3$.
For the cases $k\in\{1, 2\}$, $\Big\lfloor{\frac{k+1}{4}}\Big\rfloor= 0$, so $\sin(\frac{2\pi}{n}(m+\lfloor{\frac{k+1}{4}}\rfloor)) = \sin\Big(\frac{2\pi}{n}m\Big)$, but $m=\frac{(n-k)}{4}$, so 
$\sin(\frac{2\pi}{n}m)= \sin(\frac{2\pi(n-k)}{4n})=\sin(\frac{\pi}{2} - \frac{2\pi k}{4n}) = \cos(\frac{\pi k}{2n})$.

For the case $k=3$, $\Big\lfloor{\frac{k+1}{4}}\Big\rfloor= 1$, so following a similar reasoning, 
$\sin(\frac{2\pi}{n}(m+\lfloor{\frac{k+1}{4}}\rfloor)) = \cos(\frac{\pi}{2n})$.

\begin{equation}\label{eq:tan}
\pi_n = \begin{cases} 
          n \tan\left(\frac{\pi}{n}\right) & n=4m \\
          n \tan\left(\frac{\pi}{n}\right)\cos(\frac{\pi}{2n}) & n=4m+1,\, n=4m+3 \\
          n \tan\left(\frac{\pi}{n}\right)\cos(\frac{\pi}{n}) =n\sin\left(\frac{\pi}{n}\right) & n=4m+2.
\end{cases}
\end{equation}
Note that the condition for the second case can be written as $n=2m+1$.

\begin{figure}[htb]
	\centering
	\includegraphics[width=0.7\textwidth]{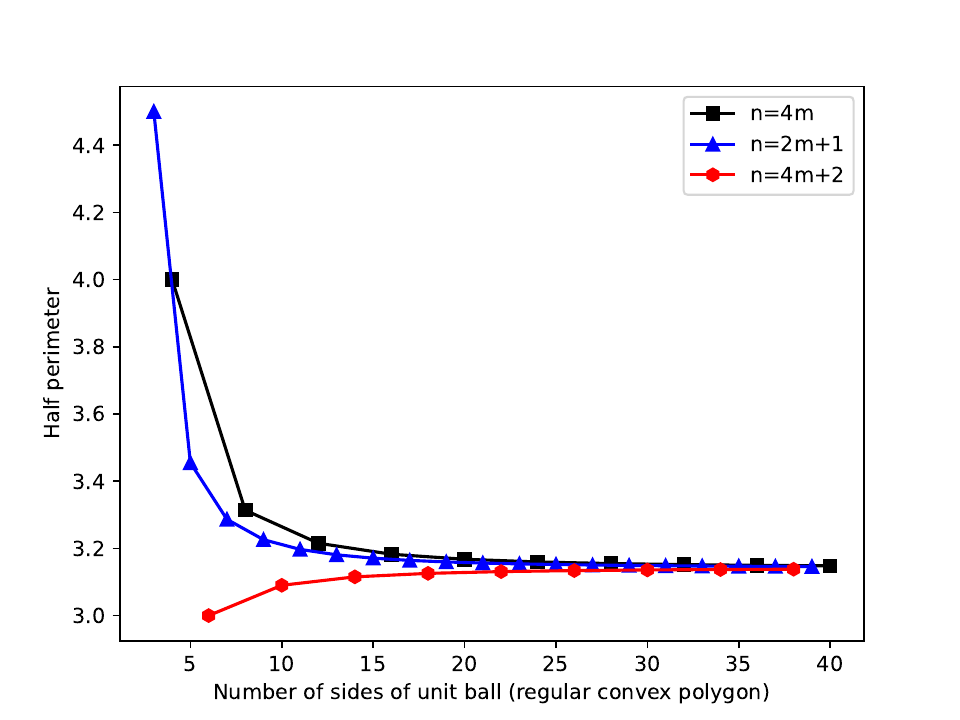}
	\caption{Values of $\pi_{B_n}$ vs. number of sides\label{fig:piscol}}
\end{figure}

Seeing that the value of $\pi_n$ is a piecewise function, we can separate the plot in Figure~\ref{fig:pis}
for each of the cases in~\eqref{eq:tan} as can be see in Figure~\ref{fig:piscol}. This clarifies the plot and gives us the intuition to
calculate bounds later.
\subsection{Beraha numbers}

The $n$th Beraha number or Beraha constant~\cite{beraha} is given by 
$B_n=2+2\cos\Big(\frac{2\pi}{n}\Big) = 4\cos\Big(\frac{\pi}{n}\Big)^2$.

The sequence\footnote{$\phi$ is the golden ratio, $\phi=\frac{1+\sqrt{5}}{2}$ and
$S$ is the silver constant.} is shown in Table~\ref{tab:ber},

\begin{table}[hbt]
\begin{tabular}{ll}
n&  $B_n$ \\
\hline
$1$&$4$\\
$2$&$0$\\
$3$&$1$\\
$4$&$2$\\
$5$&$\phi+1$\\
$6$&$3$\\
$7$& $2+2\cos(\frac{2\pi}{7})=S$\\
$8$&$2 + \sqrt{2}$\\
$9$&$2+2\cos(\frac{2\pi}{9})$\\
$10$&$\phi+ 2$,\\
\end{tabular}
	\caption{Beraha numbers $B_n$ \label{tab:ber}}
\end{table}

There are simple algebraic expressions of $\pi_n$ in terms of Beraha constants.

When $n=4m$, applying the
pythagorean identity:
$\pi_n = n\tan\Big(\frac{\pi}{n}\Big) = n \frac{\sin\Big(\frac{\pi}{n}\Big)}{\cos\Big(\frac{\pi}{n}\Big)} =  n \frac{\sqrt{1-\cos\Big(\frac{\pi}{n}\Big)^2}}{\cos\Big(\frac{\pi}{n}\Big)}$.
Applying the half angle identity, $\cos(a/2)=\sqrt{\frac{1+\cos(2a)}{a}}$ and simplifying,
$\pi_n = n\sqrt{\frac{1-\cos(\frac{2\pi}{n})}{1+ \cos(\frac{2\pi}{n})}} = n\sqrt{\frac{4-B_n}{B_n}}$.

Applying the half angle identity to the extra terms for the other cases:
$\cos\Big(\frac{\pi}{n}\Big) = \sqrt{\frac{1+\cos\Big(\frac{2\pi}{n}\Big)}{2}} =\sqrt{\frac{B_n}{4}}$
and
$\cos\Big(\frac{\pi}{2n}\Big) = \sqrt{\frac{1+\cos\Big(\frac{\pi}{n}\Big)}{2}} = \sqrt{\frac{1+\sqrt{\frac{B_n}{4}}}{2}}$.

Putting it all together:

\begin{equation*}
\pi_n = \begin{cases} 
          n\sqrt{\frac{4-B_n}{B_n}}& n=4m \\
           n\sqrt{\frac{4-B_n}{4B_n}(2+\sqrt{B_n})} =  n\sqrt{\frac{B_{2n}(4-B_n)}{4B_n}} & n=4m+1,\, n=4m+3 \\
          n\sqrt{\frac{4-B_n}{4}}& n=4m+2.
\end{cases}
\end{equation*}

\subsection{Polygonal circle numbers}
In~\cite{circlenumbers} so-called polygonal circle numbers were defined. They correspond to a different
generalization of $\pi$ from the one defined in this paper. They define circle numbers as
$c_n = \alpha_B (B) = \int_{B} ||x||_B $, the area of the unit ball calculated radially pointing outside. They find them to be $c_n = \pi \frac{\sin\Big(\frac{2\pi}{n}\Big)}{\frac{2\pi}{n}}$.

Simplifying and applying the double integer identity, $\sin(2a)=2\sin(a)\cos(a)$,
$c_n = \frac{n \sin(\frac{2\pi}{n})}{2} = n \sin\Big(\frac{\pi}{n}\Big)\cos\Big(\frac{\pi}{n}\Big)$.

We can write $\pi_n$ as a function of $c_n$,
\begin{equation}
\pi_n = \begin{cases} 
          c_n \frac{1}{(\cos\left(\frac{\pi}{n}\right))^2} & n=4m \\
          c_n \frac{\cos\left(\frac{\pi}{2n}\right)}{(\cos\left(\frac{\pi}{n}\right))^2}  & n=4m+1,\, n=4m+3 \\
          c_n \frac{1}{\cos\left(\frac{\pi}{n}\right)}  & n=4m+2.
\end{cases}
\end{equation}
\subsection{Apothem}
The radius of the inscribed circle or apothem~\cite{coxeter} of a polygon of circumscribed circle of radius one
is $r_n = 2\sin(\frac{\pi}{n})$.

\begin{equation}
\pi_n = \begin{cases} 
          \frac{n r_n}{2\cos\left(\frac{\pi}{n}\right)} & n=4m \\
          \frac{n r_n \cos(\frac{\pi}{2n}) }{2\cos\left(\frac{\pi}{n}\right)} & n=4m+1,\, n=4m+3 \\
          \frac{n r_n}{2} & n=4m+2.
\end{cases}
\end{equation}
\subsection{Bounds for $\pi_n$}
\begin{thm}
The value of $\pi_n$ is bounded:
\begin{equation}
\pi_n \in \begin{cases} 
          (\pi,4] & \text{when }n=4m \\
          (\pi,\frac{9}{2})  & \text{when }n=4m+1,\, n=4m+3 \\
            [3,\pi)   & \text{when }n=4m+2.
\end{cases}
\end{equation}

\end{thm}
\begin{proof}
We can find the maximum and minimum for $\pi_n$ by analyzing equation~\eqref{eq:tan}.

We can start with the case $4m+2$. This case,  $f(n) = n\sin\left(\frac{\pi}{n}\right)$, grows uniformly
when $n\ge 3$. With the Taylor series, it is easy to see that the derivative is positive
and the function increases monotonically.
Taking into account $\sin\left(\frac{\pi}{n}\right)  = \frac{\pi}{n} + O(1/n^3)$, and that $\pi_n=n\sin\left(\frac{\pi}{n}\right)$ then, the $\lim_{n \to \infty}\pi_n =\pi$.
As it grows monotonically with $n$, it attains
its mimimum  on the value $\pi_6=3$, so $3 \le \pi_{4m+2} \le \pi$.

The case $n=4m$,  $g(n) = n \tan\left(\frac{\pi}{n}\right)$ decreases monotonically when $n>=4$,
because $n$ increases but $\tan\left(\frac{\pi}{n}\right)$ decreases and the second term
dominates as can be seen from the Taylor series.
All the terms are positive, so for the derivative, all the terms are negative when $n>0$.

The limit $\lim_{n \to \infty}\pi_n =\pi$ is the same as the last case divided by $\cos\left(\frac{\pi}{n}\right)$ which tends to $1$, so the limit is still $\pi$.
Because of this, it attains it maximum when $\pi_4=4$ so $4 \le \pi_{4m} \le \pi$.

Finally, the cases where n is odd also decrease uniformly. It can
be seen as the case for $n=4m+2$ but with a faster cosine, which is the term which dominated.
The limit is still $\pi$ because $\lim_{n \to \infty}\cos\left(\frac{\pi}{2n}\right)=1$, so $\frac{9}{2} \ge \pi_{2k+1} \ge \pi$, where $k \ge 1$, with its
maximum at $\pi_3=\frac{9}{2}$.

\begin{equation}
\pi_n \in \begin{cases} 
          (\pi,4] & \text{when }n=4m \\
          (\pi,\frac{9}{2})  & \text{when }n=4m+1,\, n=4m+3 \\
            [3,\pi)   & \text{when }n=4m+2.
\end{cases}
\end{equation}
\end{proof}
Note that the first case, when the polygons have a number of sides $n=4m$ defines a (symmetric) norm with quarter-turn symmetry, which is known to have
$\pi_{4m} \ge \pi$ as is proven in~\cite{extremal}.

The other family of convex regular polygons inducing (symmetric) norms,
are the ones with number of sides $n=4m+2$. As stated by the following theorem these norms belong
to $\mathcal{R}$ (the set of norms coming from a linear isomorphism of a Radon norm),
which by Theorem~\ref{radonpi} have $\pi_{4m+2} \in [3, \pi)$, and they are the only ones in that range, so no other
convex regular polygons can induce a norm in $\mathcal{R}$.

\begin{thm}
Norms with convex regular polygons of $n=4m+2$ sides as unit ball are Radon.
\end{thm}
\begin{proof}
To prove that these polygons induce Radon norms, 
\begin{figure}[htb]
	\centering
	\includegraphics[width=0.5\textwidth]{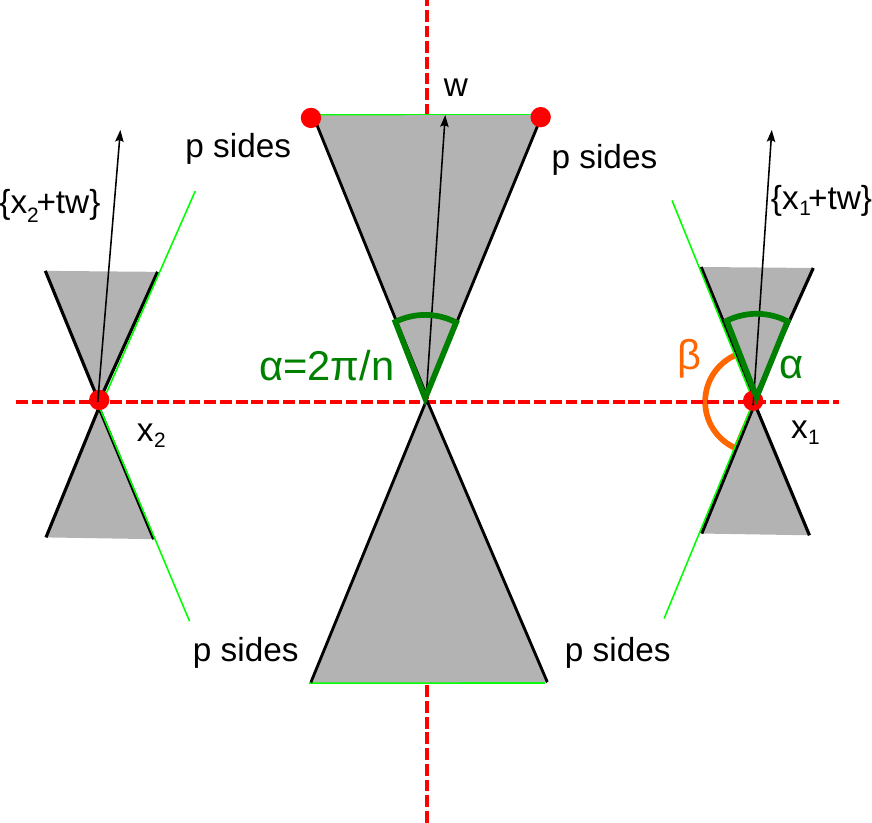}
	\caption{Polygon inducing radon norm\label{fig:radonpolyg}}
\end{figure}
we can follow a similar proof than Example 4.6 in~\cite{extremal}.
By Definition~\ref{radon}, we need to prove
Birkhoff orthogonality to be symmetric.  We only need to prove it for the
boundary of the unit ball, $\partial B$. If we take a point $w$ on an edge (or its opposite), it admits only
one supporting line (one line $u+tw$ intersecting $\partial B$ such that $||w+tu||_B \ge ||w||_B$): the side to which it belongs. Translating this
line to the origin, we find two points, vertices, $x_1$ and $x_2$ for
which $ x_1 \vdash_X w$ and $ x_2 \vdash_X w$. These vertices admit many
supporting lines (of the kind $||x_n+tw||_B \ge ||x_n||_B$) depicted as
a grey wedge in Figure~\ref{fig:radonpolyg}. They intersect the $4$ sides adjacent to
the $2$ vertices.
The angle $\alpha$ depicted in the figure, is $\alpha=\frac{2\pi}{n}$ because the
internal angle of a polygon is $\beta=\frac{(n-2)\pi}{n}$, so the angle of the wedges
is $\alpha=\pi-\beta=\pi-\frac{(n-2)\pi}{n} =\frac{2\pi}{n}$.

Each of these two wedges, when
translated to the origin, intersect with the original edge and no other, so  $w\vdash_X
x$. In other words, we have  each $4=2+2$ sides symmetrically Birkhoff orthogonal
to the corresponding $2=1+1$ and $4p$, with $p\ge0$, sides (the rest of the sides)
which are not Birkhoff orthogonal in any direction with them. This makes
in total $4+2+4p = 4m+2$ sides.  Another proof of this can be found
in~\cite{mandal2019geometric}.
\end{proof}

The three families of polygons are then:
\begin{equation}
\pi_n \in \begin{cases} 
          (\pi,4] & \text{when }n=4m:\, \text{Quarter-turn symmetry,} \\
          (\pi,\frac{9}{2})  & \text{when }n=4m+1,\, n=4m+3:\, \text{Asymmetric norms,}  \\
            [3,\pi)   & \text{when }n=4m+2:\, \text{Radon norms} .
\end{cases}
\end{equation}

\begin{corollary}
Combining all three families of polygons, $\pi_n \in [3,\frac{9}{2}]$.
\end{corollary}

\section{Polygons with an offset center\label{sec:offset}}

In the following section we analyse the behaviour of $\pi_B$ for some convex polygons with an
offset center. First we consider the general case for a triangle (which has to be isosceles).
Then, we analyse the cases for regular squares and hexagons with an offset center.

\subsection{Offset isosceles triangle}

Let's relax the conditions as much as possible and study the half-perimeter of an isosceles triangle
using the offset Minkowski functional from definition~\ref{def:minkcent}. We will denote
it $\pi_t$. 

We refer the interested reader to the Appendix~\ref{sec:isos} for the calculations. Given an offset $h_2$ from the center,
the formula for $\pi_t$ is:

\begin{equation}\label{eq:tikiX}
\pi_t = \frac{1}{2} \Bigg(\frac{2h}{h_2}+ \frac{h}{h-h_2} + \frac{2h}{h_2}\Bigg).
\end{equation}

The derivative with respect to $h_2$ is:
$f'(h_2)= \frac{h}{2(h-h_2)^2} - \frac{2h}{h_2^2}$.

When $\frac{2h}{3} < h_2 < h$, $f'(h_2) > 0$, so $\pi_t$ increases monotonically with $h_2$. When $0 < h_2 < \frac{2h}{3} $, $f'(h_2) < 0$, so $\pi_t$ decreases monotonically with $h_2$.

Also, when $h_2 \to h$, $\pi_t  \to \infty$ and, as we saw before,
when $h_2 \to 0$, $\pi_t  \to \infty$.

\begin{figure}[htb]
	\centering
	\includegraphics[width=0.3\textwidth]{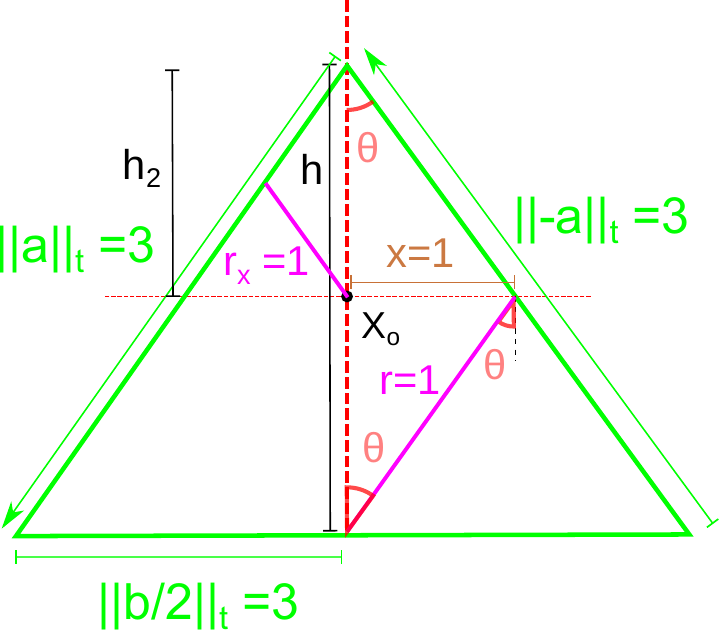}
	\caption{Offset isosceles triangle with $\pi_t=4.5$\label{fig:isos3}}
\end{figure}

The minimum value for $\pi_t$  then happens when the derivative is $0$ in the interval $0 < h_2 < h$, and can be obtained by substituting
$h_2=\frac{2h}{3}$ in equation~\eqref{eq:tiki},
\begin{equation}
\pi_t =\frac{1}{2} \Bigg(\frac{2h}{h_2}+\frac{h}{h-h_2} + \frac{2h}{h_2}\Bigg)=\frac{1}{2} \Bigg(3+3+3\Bigg) = \frac{9}{2},
\end{equation}
and
$\frac{a}{r_x} =3$,
$\frac{a}{r} =3$,
and
$\frac{b}{x} =3$.
Note that the minimum value does not depend on $a$, $b$ or any function of both.
An example of this case is depicted for a concrete isosceles triangle in Figure~\ref{fig:isos3}.

We can also calculate, for example, when $\frac{a}{r_x}<\frac{a}{r}$:
$\overset{max}{\mu}_t(B_t)  = \frac{1}{2} \Big(\frac{a}{r} + \frac{a}{r} + \frac{b}{x}\Big)$
and
$\overset{min}{\mu}_t(B_t)  = \frac{1}{2} \Big(\frac{a}{r_x} + \frac{a}{r_x} + \frac{b}{x}\Big)$.

In this case, if we set  $h_2=\frac{4h}{5}$,
$||-a||_t = \frac{a}{r_x} =\frac{2h}{h_2} = \frac{5}{2}$,
$||a||_t = \frac{a}{r} =\frac{h}{h-h_2} = 5$,
and
$||b||_t =\frac{b}{x} =\frac{2h}{h_2} = \frac{5}{2}$,
so
$\overset{max}{\mu}_t(B_t)  = \frac{25}{2}$,
$\overset{min}{\mu}_t(B_t)  = \frac{15}{2}$,
and
$\mu_t(B_t)= 10$,
so

$\pi_t = \frac{{\mu}_t(B_t)}{2}  = \frac{10}{2} = 5$.

Note that both lengths of the side $a$, $||\vec{a}||_t$ and $||-\vec{a}||_t$ are not the same, as was the case
for the centered regular triangle. Moreover, $\overset{max}{\mu}_t(B_t) \neq \overset{min}{\mu}_t(B_t)$.

Watching Figure~\ref{fig:isos3}, a general symmetry argument can be made: the minimum
happens only when $||\vec{a}||_t$ = $||-\vec{a}||_t$ and
$\overset{max}{\mu}_t(P) = \overset{min}{\mu}_t(P)$, as was the case for the centered regular triangle.

We capture the results of this section in the following observation:
\begin{obs}
Equation~\eqref{eq:tikiX} is symmetric, whether the offset $h_2$ goes up or down (see Appendix~\ref{sec:isos} for more details)
the equation is valid.
The consequence is that for every value of
$\pi_t > 4.5$, every isosceles triangle will have two corresponding center positions, with the same offset with respect to
the center ( placed in the axis of symmetry, one towards the vertex
in the axis of symmetry and one towards its opposite side $d$ and $-d$).
\end{obs}

\begin{figure}[htb]
	\centering
	\includegraphics[width=0.15\textwidth]{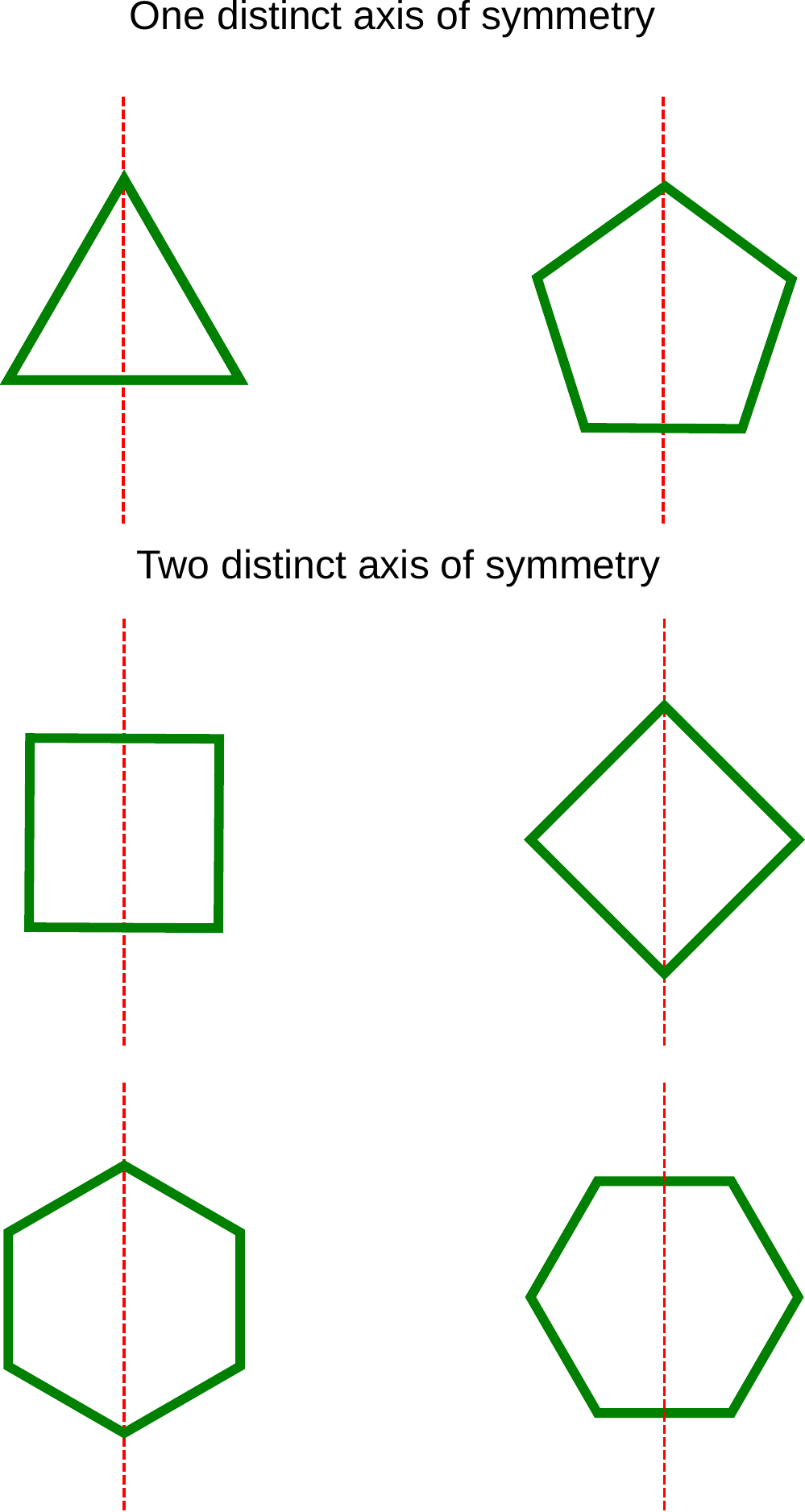}
	\caption{Number of non-equivalent axes of symmetry\label{fig:axis}}
\end{figure}

For polygons with an even number of vertices, there are two possible
non equivalent axes of symmetry, one which crosses two vertices and
one which goes from the center of one side to the 
center of the other side. This is depicted in Figure~\ref{fig:axis}.

\subsection{Square with offset center}

\begin{figure}[htb]
	\centering
	\includegraphics[width=0.5\textwidth]{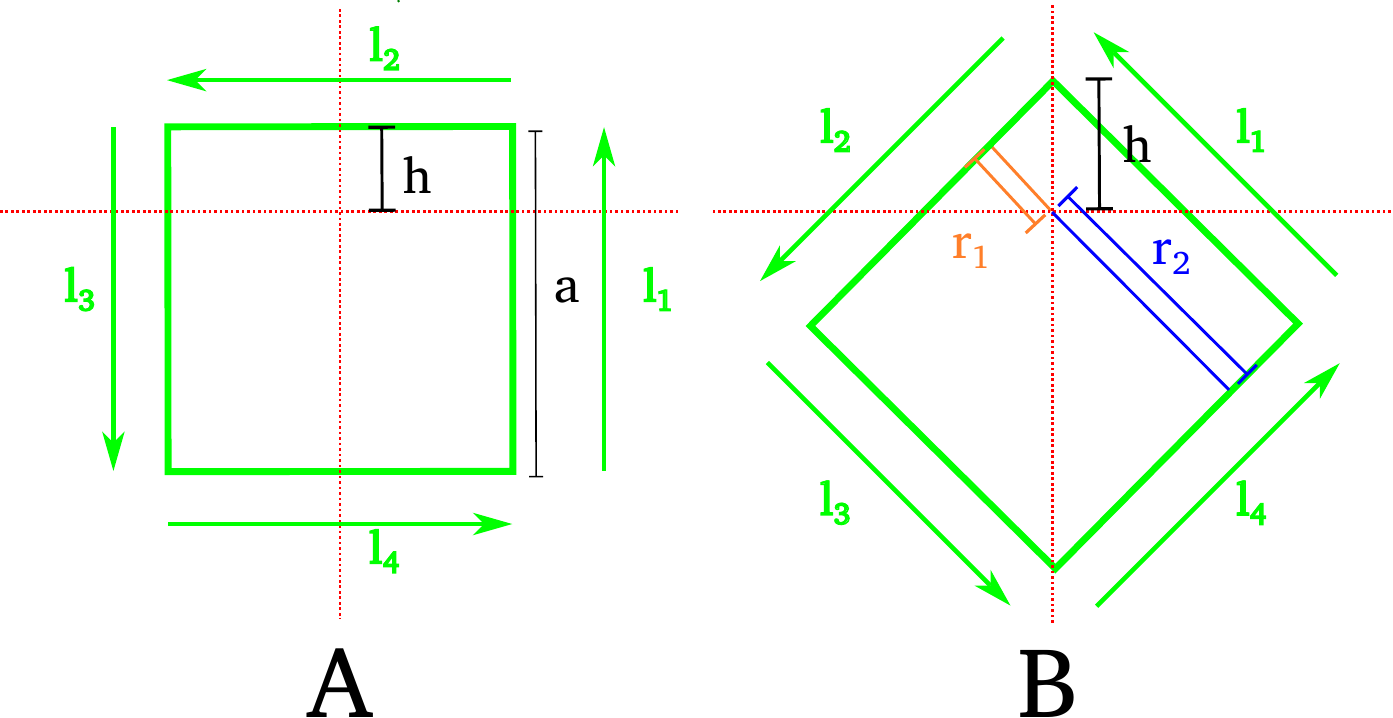}
	\caption{Offset square\label{fig:offsq}}
\end{figure}

Figure~\ref{fig:offsq} depicts the two configurations for a square.

For configuration A:
$||l_1||_A = a/h$, $||l_2||_A = 2$, $||l_3||_A = \frac{a}{a-h}$ and 
$||l_4||_A = 2$.
\begin{equation}
\pi_A = \frac{1}{2} (||l_1||_A + ||l_2||_A + ||l_3||_A + ||l_4||_A) = 2+\frac{a}{2h} + \frac{a}{2(a-h)}.
\end{equation}

Taking the derivative:
$\frac{d \pi_A }{dh} = - \frac{a}{2h^2} + \frac{a}{2(a-h)^2} $,
we can equate it to $0$ and:
$- \frac{a}{2h^2} + \frac{a}{2(a-h)^2}  = 0$, then
$h=\frac{a}{2}$.
When $h < \frac{a}{2}$, the second term dominates, because $|\frac{a}{2(a-h)^2}| <  |\frac{a}{2h^2}|$,
so the derivative is negative and $h=\frac{a}{2}$ is a minimum, which gives
\begin{equation}
min_{h=\frac{a}{2}} (\pi_A) = 2+\frac{2a}{2a} + \frac{a}{2(a-a/2)} = 4.
\end{equation}

For configuration B:
By similarity of triangles, $r_1 = \frac{h}{\sqrt{2}}$,
so
$||l_1||_B = \frac{a}{r_1} =\frac{a\sqrt{2}}{h} = ||l_2||_B$
and
$||l_3||_B = \frac{a}{a-r_1} = \frac{a}{a-\frac{h}{\sqrt{2}}}= ||l_4||_B$,

Finally
\begin{equation}
\pi_B = \frac{1}{2} (||l_1||_B + ||l_2||_B + ||l_3||_B + ||l_4||_B) = \frac{a\sqrt{2}}{h} + \frac{a}{a-\frac{h}{\sqrt{2}}}.
\end{equation}

The derivative is:
$\frac{d \pi_B }{dh} = -\frac{a\sqrt{2}}{h^2} + \frac{a}{\sqrt{2}(a-\frac{h}{\sqrt{2}})^2} $.
This derivative is the same as the one for configuration A but with $h_b = \frac{h}{\sqrt{2}}$,
so following the same reasoning, the minimum happens when $h_b = \frac{a}{2} = \frac{h}{\sqrt{2}}$, so
$h = \frac{\sqrt{2}a}{2}= \frac{a}{\sqrt{2}}$.
Finally,
\begin{equation}
min_{h=\frac{a}{\sqrt{2}}} (\pi_B) = 4.
\end{equation}

\begin{obs}
Given a value of $\pi_s > 4$, we can find four distinct configurations of the square (two per distinct
axis of symmetry) with that value of $\pi_s$ as half-perimeter.
\end{obs}
\subsection{Regular hexagon with offset center}

\begin{figure}[htb]
	\centering
	\includegraphics[width=0.8\textwidth]{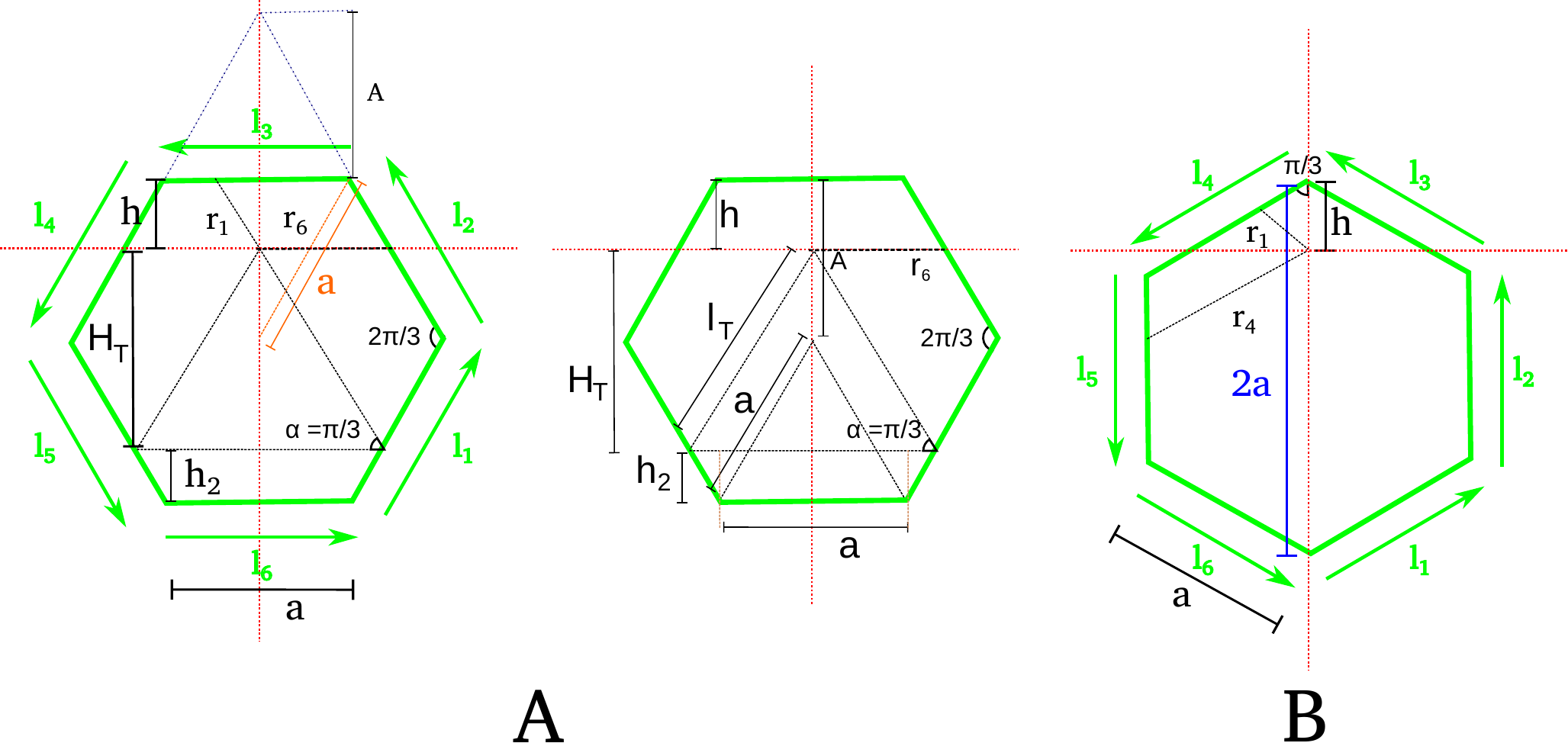}
	\caption{Offset hexagon (detail of A)\label{fig:offhex}}
\end{figure}
For the A configuration in Figure~\ref{fig:offhex}: 
The apothem of a regular hexagon, is $A=a\frac{\sqrt{3}}{2}$.
We can calculate $r_1$ by similarity of triangles using Figure~\ref{fig:offhex}:
$r_1 = a\frac{h}{A} =2\frac{h}{\sqrt{3}}$.

From triangle similarity of the top triangle,
$\frac{A+h}{A} = \frac{r_6}{a/2}=\frac{a\frac{\sqrt{3}}{2}+h}{a\frac{\sqrt{3}}{2}}$,
so $r_6=\frac{a}{2} + \frac{h}{\sqrt{3}}$, $h_2=2A-h-H_T$, and $H_T=l_t \frac{\sqrt{3}}{2}$ where
$l_T$ is the side of the (equilateral) triangle in the figure.
Substituting $H_T$ and $A$ in $h_2$, $h_2=a\sqrt{3}-h-l_T\frac{\sqrt{3}}{2}$.
From that figure $l_T = 2h_2\frac{cos(\pi/3)}{cos(\pi/6)}+a = \frac{2h_2}{\sqrt{2}} + a$. Substituting 
$h_2$, $l_T =3a-\frac{2h}{\sqrt{3}}-l_T$ so, $l_T=\frac{1}{2}(3a-\frac{2h}{\sqrt{3}})$.
Now we can calculate the length of the sides:
$||l_1||_A = \frac{a}{r_1} = \frac{a\sqrt{3}}{2h} = ||l_2||_A$,
$||l_3||_A  = ||l_6||_A  = \frac{a}{r_6} = \frac{a}{2\frac{h}{\sqrt{3}}}$
and $||l_5||_A  = ||l_4||_A  = \frac{a}{l_T} = \frac{2a}{3a-\frac{2h}{\sqrt{3}}}$.

From these, we obtain:
\begin{equation}
\pi_A = \frac{a\sqrt{3}}{2h} + \frac{a}{2\frac{h}{\sqrt{3}}} +\frac{2a}{3a-\frac{2h}{\sqrt{3}}} = \frac{a\sqrt{3}}{h} +\frac{2a}{3a-\frac{2h}{\sqrt{3}}} .
\end{equation}

To find the minimum, we will simplify it with a substitution $H=\frac{2h}{\sqrt{3}}$.
\begin{equation}
\pi_A = \frac{2a}{H} +\frac{2a}{3a-H}.
\end{equation}
With this substitution, we find the derivative:
\begin{equation}
\frac{d \pi_A}{dH}= \frac{-2a}{H^2} +\frac{2a}{(3a-H)^2}.
\end{equation}

The formula for $\pi_A$ is invalid when $h>\frac{a\sqrt{3}}{2}$, or $H>a$,
so we cannot find the minimum by equating the derivative to $0$, but
for $0\le H\le a$ the second term is dominated by the first $|\frac{2a}{H^2}| \ge |\frac{2a}{(3a-H)^2}|$
so the derivative is negative as $H$ grows, so the minimum is:

\begin{equation}
min_{B = a}(\pi_A)  = min_{h=\frac{a\sqrt{3}}{2}} (\pi_A) = 3.
\end{equation}

For the B configuration:
Seeing the figure,
$r_1 cos(\frac{\pi}{6}) = \frac{h}{2}$, so $r_1 = h$.

Then
$||l_1||_B = \frac{a}{r_1} = \frac{a}{h} $,
$||l_2||_B = \frac{a}{h} $,
$||l_3||_B = ||l_1||_B =  \frac{a}{h} $,
$||l_4||_B = 1$,
$||l_5||_B = \frac{a}{2a-h} $ and
$||l_6||_B = ||l_1||_B =||l_4||_B =  1$,
so we obtain:
\begin{equation}
\pi_B = \frac{1}{2} (||l_1||_B + ||l_2||_B + ||l_3||_B + ||l_4||_B) = 1+\frac{a}{2a-h} + \frac{a}{h}.
\end{equation}
To find the minimum:
$\frac{d \pi_B }{dh} = \frac{a}{(2a-\frac{h})^2} -\frac{a}{h^2}  = 0$, so $h=a$.
We know that $|\frac{a}{(2a-h)^2}| < |\frac{a}{h}|$ when $h < a $ so the derivative is negative and we
have found a minimum so:
\begin{equation}
min_{h=a} (\pi_B) = 3.
\end{equation}

\begin{obs}\label{hexaobs}
Given a value of $\pi_h > 3$, we can find four configurations of the hexagon (two per distinct
axis of symmetry), with different positions of the center, with that value of $\pi_h$ as half-perimeter.
\end{obs}


\section{Minimum $\pi_B$ for the general case\label{sec:bounds}}

In this section we find the general bounds for $\pi_B$ for the general case of an asymmetric norm
with one axis of symmetry.

First we need the following lemma and some definitions.

\begin{lem}\label{extreme}
Given $d(y) : \R \rightarrow \R$, the width of the unit ball as a function of the displacement $y$ along the axis of symmetry.
It will have either a unique point as maximum or an interval.
The unique maximum (be it a point or an interval) may be at either ends ($d(y_{max})$ or $d(y_{min})$).
If the maximum is not at one end, that end has to be a minimum.
\end{lem}
\begin{proof}
The convexity of the unit ball and its axis of symmetry, limits the number of maxima of $d(y)$ to either
a unique point or an interval. If there where two separate maxima in $d(y)$, $d(y_1)$ and $d(y_2)$
there would be a minimum between them or $d(y_1) = d(y_2)$, by the extreme value theorem.
If there was a minimum between them,
then the segment connecting the two maxima would not be completely contained in the unit ball, so it would not be convex.
\end{proof}

From~\cite{extremal}, we borrow the following definition:
\begin{definition}
Let $B \subset \R^2$ be a convex set. If $F\subset B$ is non-empty, convex, and
$x, y \in B\textnormal{ and } (x,y) \cap F \neq 0 \Rightarrow x, y \in F$,
we call $F$ a \emph{face} of B. We say $p\in B$ is an \emph{extreme point} if $\{p\}$ is a face, i.e.
$x, y \in B\textnormal{ and } p\in (x,y) \Rightarrow x=y=p$.
Also from~\cite{extremal} we state (fact (a) after the definition)
\begin{obs}\label{proper}
If $F \subset B$ is a proper face of $B$ (i.e. $F\ne B$), then $F \subset \partial B$.
\end{obs}

\begin{obs}\label{ends}
By Lemma~\ref{extreme}, if we follow the axis of symmetry of a convex set, the two ends are either a face or an extreme point.
\end{obs}

\end{definition}

\begin{figure}[htb]
	\centering
	\includegraphics[width=0.4\textwidth]{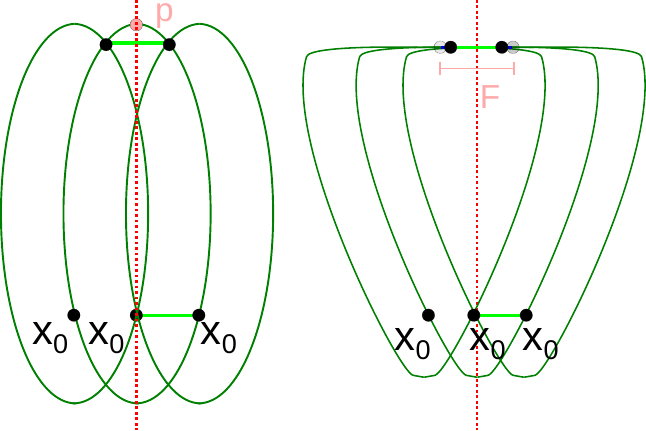}
	\caption{Two examples of the construction\label{fig:pregolab}}
\end{figure}

Finally, we are ready to prove the main theorem:

\begin{thm}
Given a convex set $B$ with one axis of symmetry, the asymmetric norm defined by~\eqref{def:minkcent} will measure its
half-perimeter $\pi_B$ to be $\pi_B \ge 3$ and $\pi_B$ can take any value in that range, depending on the concrete norm.
\end{thm}

\begin{proof}
In order to find the minimum $\pi_B$, we are going to use the construction shown in Figure~\ref{fig:mingolab}.
In that figure $x_0$ is the center. We displace the unit ball laterally until it touches the center.
At either ends following the axis of symmetry, by~\ref{ends} we either have that the end is a face, which can have
$\mu(F) \le d_{x_0}/2$ or an extreme point (it is the same case $\mu(F) = 0$) or a face with $\mu(F) > d_{x_0}/2$.
For these two cases we have the two behaviours shown in Figure~\ref{fig:pregolab}. In
the first case $\mu(F) \le d_{x_0}/2$, where the unit ball cuts its copies at the
endpoints of a segments of length $d_{x_0}/2$. If not, the length of the segment will be bigger than $d_{x_0}/2$,
but the segment will be completely contained in $\partial B$ (because of Observation~\ref{proper}), so we can take a subinterval
measuring $d_{x_0}/2$ and it will still be contained in $\partial B$.

\begin{figure}[htb]
	\centering
	\includegraphics[width=0.4\textwidth]{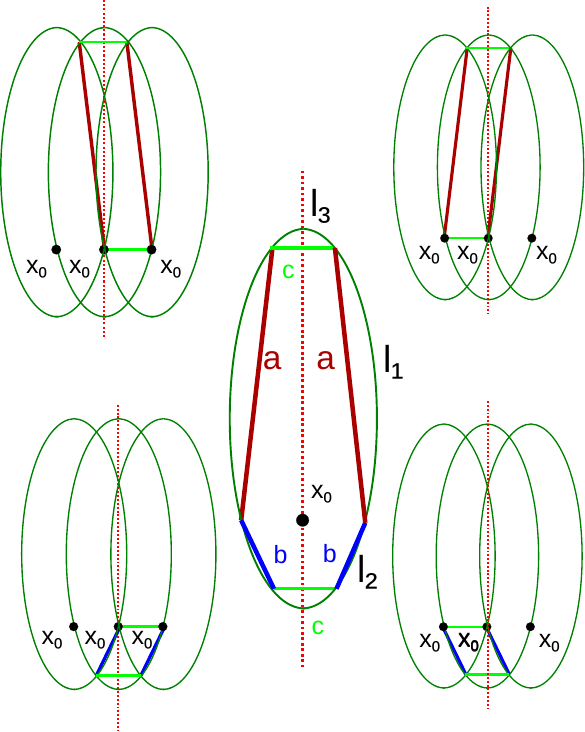}
	\caption{Inscribe an hexagon in the unit circle\label{fig:mingolab}}
\end{figure}

\begin{figure}[htb]
	\centering
	\includegraphics[width=0.2\textwidth]{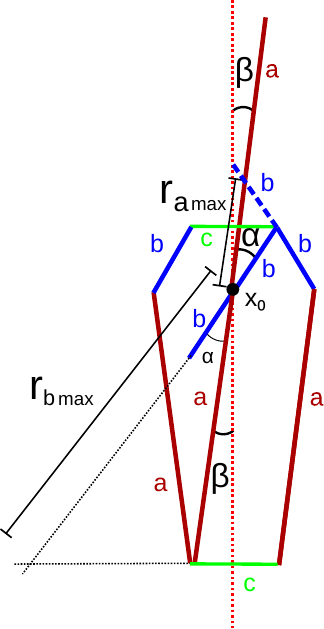}
	\caption{Constructions to find $r_{bmax}$ and $r_{amax}$ \label{fig:mingolab3}}
\end{figure}

Using the construction depicted in~\ref{fig:mingolab}, we can inscribe an hexagon $H$.

Its sides have Euclidean length $||l_1||=a$, $||l_2||=b$, $||l_3|| = c$ and  $||-l_3|| = c$,
with $c$, half the diameter at $x_0$, i.e. $c=d(x_0)/2$.
The length given by the Minkowski norm of each side is $||l_1||_B=1$, $||l_2||_B=1$, $||l_3||_B=1$ and $||-l_3||_B=1$.
Note that $||l_1||_B=1$ and $||l_2||_B=1$ because the side is $d(x_0)/2$ long and then  $||l_3||_B=1$ (or conversely
$||-l_3||_B=1$), it is a case of parallel transport as depicted in Figure~\ref{fig:mingolab}.

From Figure~\ref{fig:mingolab3}, we can find the maximum radius of the unit ball for measuring $||-l_1||_B$ and $||-l_2||_B$.
Without loss of generality, we assume $a\ge b$. Then, $r_{bmax} = \frac{acos(\beta)}{cos(\alpha+\beta)}$ and
$r_{amax}=\frac{bcos(\alpha+\beta)}{cos(\beta)}$. From this:
$||-l_1||_B \le  \frac{a}{r_{amax}} = \frac{a\cos(\beta)}{b\cos(\alpha+\beta)}$,
\mbox{$||-l_2||_B \le \frac{b}{r_{bmax}}=\frac{b\cos(\alpha+\beta)}{a\cos(\beta)}$}.
We can name  $\frac{a\cos(\beta)}{b\cos(\alpha+\beta)} = y$.
with $y>0$. The half-perimeter of the hexagon is then:
$\frac{\mu_B(H)}{2} =\frac{||l_1||_B+||-l_1||_B+ ||l_2|_B+||-l_2||_B+ ||l_3||_B+||-l_3||_B}{2}\ge$
$\frac{4+y+ \frac{1}{y}}{2}$,
so $\frac{\mu_B(H)}{2}  \ge 2+\frac{y+ \frac{1}{y}}{2} = f(y)$,
This function has an extremum at, $f'(y) = \frac{1- \frac{1}{y^2}}{2}=0$, so, $y=1$, and $\frac{\mu_B(H)}{2}  = 3$.
It is a minimum, because $f''(y)  = \frac{1}{y^3}>0$.
From the triangle inequality we have proposition~\eqref{inscribed}, so we know that $\pi_B=\frac{1}{2}\mu_B(\partial B)\ge \frac{\mu_B(H)}{2} \ge 3 $ and we can conclude that for any asymmetric Minkowski space with an axis of symmetry,
the minimum value  of $\pi_B$, $min(\pi_B) = 3$.
From Observation~\ref{hexaobs} we can conclude that $\pi_B$ takes all the possible values in the range $\pi_B \ge \pi$.
\end{proof}

\section*{Acknowledgements}

Thanks to Germain Fruteau for may fruitful conversations and computer simulations which have helped this paper immensely.
\bibliographystyle{vancouver}
\bibliography{pi}
%
%
%

\appendix
\section{Appendix}
\subsection{Calculations for the regular pentagon \label{sec:penta}}

For the regular pentagon, the apothem $a_5=\frac{l}{2\tan(\frac{\pi}{5})}$, and $||l||_{\pentagon} = 2-2 a_5 \tan(\frac{\pi}{10})=$
$2-\frac{||l||_{\pentagon} }{\tan(\frac{\pi}{5})} \tan(\frac{\pi}{10})$ so $||l||_{\pentagon} =\
\frac{2}{1+\frac{\tan(\frac{\pi}{10})}{\tan(\frac{\pi}{5})}} =\
 \frac{(5-\sqrt{5})}{2}$. Finally $\pi_{\pentagon} = \frac{5||l||_{\pentagon}}{2} = \frac{5(5-\sqrt{5})}{4}$.

\subsection{Calculations for the offset isosceles triangle \label{sec:isos}}
\begin{figure}[htb]
	\centering
	\includegraphics[width=0.5\textwidth]{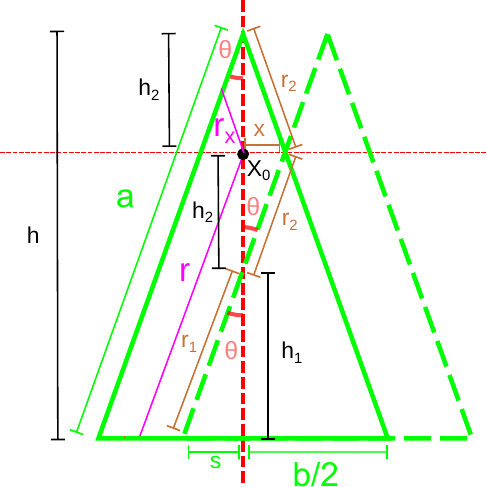}
	\caption{Offset isosceles triangle, first configuration\label{fig:isos}}
\end{figure}

Note that $\pi_t$ is actually a function of the length of the sides and the offset of the center,
$\pi_t(a, b, x_o)$.
We refer to Figure~\ref{fig:isos} in this section. We will skip, for the sake of clarity
 in the following
equations the Euclidean norm, so when we denote $r$ we actually mean $||\vec{r}||$ and so on.
First,
\begin{equation}\label{eq:pitval}
\pi_t = \frac{1}{2} \Big(\frac{a}{r} +\frac{a}{r_x}+ \frac{b}{x}\Big).
\end{equation}
We define $h$ as
\begin{equation}\label{eq:h}
h = h_1+2h_2.
\end{equation}
and applying Pythagoras:
\begin{equation}\label{eq:hab}
h = \sqrt{a^2 - \frac{b^2}{4}}.
\end{equation}
We also have
\begin{equation}\label{eq:rx}
\frac{r_2}{2x} = \frac{r_x}{x},
\end{equation}
\begin{equation}\label{eq:two}
\frac{x}{h_2} = \tan(\theta)
\end{equation}
and
\begin{equation}\label{eq:three}
\frac{s}{h_1} = \tan(\theta),
\end{equation}
but also,
\begin{equation}\label{eq:four}
\frac{\frac{b}{2}}{h} = \frac{b}{2h} = \tan(\theta).
\end{equation}
Also from Pythagoras,
\begin{equation}\label{eq:five}
s^2+h_1^2 = r_1^2
\end{equation}
and
\begin{equation}\label{eq:six}
x^2+h_2^2=r_2^2.
\end{equation}

Note that $a$ and $b$ are the dimensions of the triangle and $h_2$ set the position of
the origin. All the other dimensions depend on these per the equations above.

Combining equations~\eqref{eq:two} and~\eqref{eq:four},
\begin{equation*}
\frac{x}{h_2} = \frac{b}{2h},
\end{equation*}
so
\begin{equation}\label{eq:xval}
x= \frac{bh_2}{2h}.
\end{equation}
Similarly, from~\eqref{eq:two} and~\eqref{eq:three},
\begin{equation}
s=\frac{xh_1}{h_2}
\end{equation}
and, from~\eqref{eq:five}
\begin{equation}\label{eq:comb}
\frac{h_1^2x^2}{h_2^2}+h_1^2=r_1^2.
\end{equation}
We also have
\begin{equation}
r=r_1 + r_2,
\end{equation}
and from~\eqref{eq:six}, combined with~\eqref{eq:comb}
\begin{equation*}
r=\sqrt{\frac{h_1^2x^2}{h_2^2}+h_1^2}+\sqrt{h_2^2+x^2}=h_1\sqrt{\frac{x^2}{h_2^2}+1}+\sqrt{h_2^2+x^2}
\end{equation*}
substituting $x$ from~\eqref{eq:xval},
\begin{equation*}
r=h_1\sqrt{\frac{x^2}{h_2^2}+1}+\sqrt{h_2^2+\frac{b^2 h_2^2}{4h^2}}=h_1\sqrt{\frac{b^2}{4h^2}+1}+h_2\sqrt{1+\frac{b^2}{4h^2}},
\end{equation*}
but, from~\eqref{eq:h}, $h=h_1+2h_2$,
\begin{equation}\label{eq:rval}
r=(h-2h_2)\sqrt{\frac{b^2}{4h^2}+1}+h_2\sqrt{1+\frac{b^2}{4h^2}} = (h-h_2)\sqrt{\frac{b^2}{4h^2}+1}.
\end{equation}
Finally,~\eqref{eq:xval} and~\eqref{eq:rval} and~\eqref{eq:pitval},
\begin{equation}\label{eq:hhh}
\pi_t = \frac{1}{2} \Bigg(\frac{a}{r_x}+\frac{a}{(h-h_2)\sqrt{\frac{b^2}{4h^2}+1}} + \frac{b}{\frac{bh_2}{2h}}\Bigg)= \frac{1}{2} \Bigg(\frac{a}{r_x}+\frac{a}{(h-h_2)\sqrt{\frac{b^2}{4h^2}+1}} + \frac{2h}{h_2}\Bigg).
\end{equation}
We can substitute $h$, from~\eqref{eq:hab},
\begin{equation*}
\pi_t = \frac{a}{2r_x}+\frac{a}{2(\sqrt{a^2 - \frac{b^2}{4}}-h_2)\sqrt{\frac{b^2}{4(a^2 - \frac{b^2}{4})}+1}} + \frac{\sqrt{a^2 - \frac{b^2}{4}}}{h_2}.
\end{equation*}

But, from~\eqref{eq:rx}, $r_x=\frac{r_2}{2}$ and
from \eqref{eq:six}
$r_2 = \sqrt{x^2-h_2^2}$, so 
\begin{equation*}
r_x=\frac{\sqrt{x^2-h_2^2}}{2}=\frac{\sqrt{\frac{(bh_2)^2}{4h^2}-h_2^2}}{2}=\frac{1}{4}\sqrt{\frac{(bh_2)^2}{h^2}-4h_2^2}=\frac{h_2}{4}\sqrt{\frac{b^2}{h^2}-4},
\end{equation*}
substituting in $r_x$ and $h$, from~\eqref{eq:hab},
\begin{equation*}
\pi_t = \frac{2a}{h_2\sqrt{\frac{b^2}{a^2 - \frac{b^2}{4}}-4}}+\frac{a}{2(\sqrt{a^2 - \frac{b^2}{4}}-h_2)\sqrt{\frac{b^2}{4(a^2 - \frac{b^2}{4})}+1}} + \frac{\sqrt{a^2 - \frac{b^2}{4}}}{h_2},
\end{equation*}
\begin{equation*}
\pi_t = \frac{1}{h_2\sqrt{\frac{1}{a^2 - \frac{b^2}{4}}}}+\frac{a}{2(\sqrt{a^2 - \frac{b^2}{4}}-h_2)\frac{a}{\sqrt{a^2-\frac{b^2}{4}}}} + \frac{\sqrt{a^2 - \frac{b^2}{4}}}{h_2},
\end{equation*}
substituting back $h$, from~\eqref{eq:hab},
\begin{equation*}
\pi_t = \frac{h}{h_2}+\frac{a}{2(h-h_2)\frac{a}{h}} + \frac{h}{h_2} ,
\end{equation*}
so we obtain
\begin{equation}\label{eq:pih2h}
\pi_t = \frac{h}{h_2}+\frac{h}{2(h-h_2)} + \frac{h}{h_2}.
\end{equation}

Note that $\pi_t$ is not bounded. When $h_2 \to 0$, $\pi_t \to \infty$.

The other end, cannot be represented with this construction, because $h_1$ would be negative.

\begin{figure}[htb]
	\centering
	\includegraphics[width=0.7\textwidth]{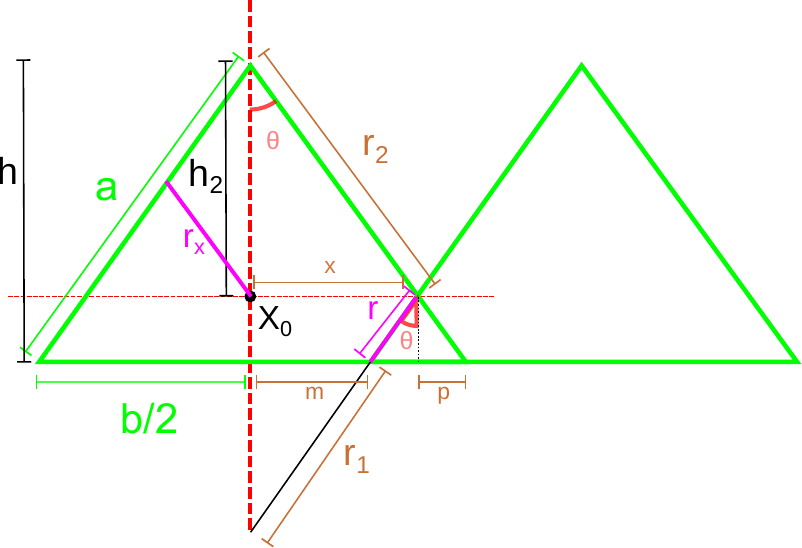}
	\caption{Offset isosceles triangle, second configuration\label{fig:isos2}}
\end{figure}

We show the second configuration, in Figure~\ref{fig:isos2}.
In this configuration, the next are true:
\begin{equation}\label{eq:eqb2}
\frac{h-h_2}{p}=\frac{h_2}{x},
\end{equation}
\begin{equation}\label{eq:rx2}
\frac{r_2}{2x} = \frac{r_x}{x},
\end{equation}

\begin{equation}\label{eq:eqd2}
r_1+r=r_2,
\end{equation}
\begin{equation}\label{eq:eqe2}
r_2=a-r,
\end{equation}
and, from Pythagoras,
\begin{equation}\label{eq:eqf2}
r^2 = p^2+(h-h_2)^2,
\end{equation}
\begin{equation}\label{eq:hab2}
h = \sqrt{a^2 - \frac{b^2}{4}},
\end{equation}
and,
\begin{equation}\label{eq:eqg2}
r_2^2 = x^2+h_2^2.
\end{equation}

From~\eqref{eq:eqf2} and~\eqref{eq:eqe2},
\begin{equation}\label{eq:simplex2}
x=\sqrt{(a-r)^2-h_2^2}
\end{equation}

Combining~\eqref{eq:eqf2},~\eqref{eq:simplex2}  and~\eqref{eq:eqb2},
\begin{equation*}
r^2=\frac{x^2(h-h_2)^2}{h_2^2}+(h-h_2)^2 = \frac{x^2(h-h_2)^2}{h_2^2}+(h-h_2)^2,
\end{equation*}
and, substituting~\eqref{eq:simplex2},
\begin{equation*}
r^2-\frac{(h-h_2)^2(a-r)^2}{h_2^2}=0,
\end{equation*}
and expanding the right side,
\begin{equation*}
r^2\Big(\frac{h_2^2}{(h-h_2)^2}-1\Big)+2ar-a^2=0.
\end{equation*}
Solving the quadratic equation for $r$,
\begin{equation*}
r=\frac{-2a\pm\sqrt{4a^2+4a^2(\frac{h_2^2}{(h-h_2)^2}-1)}}{a\Big(\frac{h_2^2}{(h-h_2)^2}-1\Big)}=\
\frac{a(-1\pm\frac{h_2}{h-h_2})}{\frac{h_2^2}{(h-h_2)^2} - 1}.
\end{equation*}
Note that $r>0$ and $h-h_2<h_2$, so $\frac{h_2^2}{(h-h_2)^2}>1$, so we can discard the
negative solution above we introduced when we squared $r$:
\begin{equation}\label{eq:r2}
r=\frac{a(-1+\frac{h_2}{h-h_2})}{\frac{h_2^2}{(h-h_2)^2} - 1} = \frac{a}{1+\frac{h_2}{h-h_2}} .
\end{equation}

\begin{equation*}
\pi_t = \frac{1}{2} \Big(\frac{a}{r_x} + \frac{a}{r} + \frac{b}{x}\Big),
\end{equation*}
substituting $r$ from~\eqref{eq:simplex2} and $x$ from~\eqref{eq:r2},

\begin{equation*}
\pi_t = \frac{1}{2} \Bigg(\frac{a}{r_x} + \frac{h}{h-h_2} + \frac{b}{h_2\sqrt{\Big(\frac{a }{h}\Big)^2-1}}\Bigg).
\end{equation*}
substituting~\eqref{eq:hab2},
\begin{equation*}
\pi_t = \frac{1}{2} \Bigg(\frac{a}{r_x} + \frac{h}{\sqrt{a^2 - \frac{b^2}{4}}-h_2} + \frac{2}{h_2\sqrt{\frac{1}{a^2 - \frac{b^2}{4}}}}\Bigg).
\end{equation*}
Combining~\eqref{eq:rx2} and ~\eqref{eq:eqe2}, we get $r_x=\frac{r_2}{2}=\frac{a-r}{2}$,
so 
\begin{equation*}
\pi_t = \frac{1}{2} \Bigg(\frac{2a}{a-\frac{a}{1+\frac{h_2}{h-h_2}}} +\frac{h}{\sqrt{a^2 - \frac{b^2}{4}}-h_2} + \frac{2}{h_2\sqrt{\frac{1}{a^2 - \frac{b^2}{4}}}}\Bigg)
\end{equation*}
substituting~\eqref{eq:hab2} back,
\begin{equation}\label{eq:tiki}
\pi_t = \frac{1}{2} \Bigg(\frac{2h}{h_2}+ \frac{h}{h-h_2} + \frac{2h}{h_2}\Bigg).
\end{equation}
Note that~\eqref{eq:tiki} is the same as~\eqref{eq:pih2h}. In both equations, $h_2$ and $h$ represent the same dimensions (the height of the triangle and the offset
from its highest point to the chosen center respectively) so we can analyze
both simultaneously using the same equation.
\end{document}